\documentstyle{amsppt}
\document

\def\ZZ{{\Bbb Z}}
\def\Gm{{\Bbb G}_m}

\def\C{{\Cal{C}}}
\def\CU{{\Cal{CU}}}

\def\CA{{\Cal{A}}}
\def\CT{{\Cal{T}}}

\def\cbM{{\breve{\Cal{M}}}}
\def\G{{\Gamma}}
\def\bc{{\bold{\chi}}}
\def\tc{{\tilde{\chi}(*)}}
\def\ra{{\rightarrow}}

\def\QQ{{\Bbb Q}}

\def\Qbar{{\overline{\Bbb Q}}}

\def\RR{{\Bbb R}}
\def\ad{{\bold A}}
\def\CC{{\Bbb C}}

\def\GalE{{Gal(\overline{E}/E)}}
\def\GalEv{{Gal(\overline{E}_v/E_v)}}
\def\GalEp{{Gal(\overline{E}'/E')}}
\def\GalEvp{{Gal(\overline{E}'_v/E'_v)}}

\baselineskip=16pt
\hfuzz=5pt

\define\isoarrow{{\overset\sim\to\longrightarrow}}
\def\-{{-1}}

\magnification=\magstep1
\NoBlackBoxes

\centerline{\bf The local Langlands conjecture for $GL(n)$ over a $p$-adic
field, $n < p$.}
\bigskip
\centerline{\bf Michael Harris\footnote{U.R.A. 748 du CNRS et
Institut de Mathematiques de Jussieu-UMR CNRS 9994.
Supported in part by the National Science
  Foundation, through Grant DMS-9423758.}}

\bigskip

\centerline{\bf Introduction}
\medskip

Let $F$ be a $p$-adic field and $n$ a positive integer.  Let 
$\Cal{A}(n,F)$ denote
the set of equivalence classes of irreducible admissible 
representations of $GL(n,F)$,
and let $\Cal{A}_0(n,F)$ be its subset of supercuspidal representations.
Similarly, let $\Cal{G}(n,F)$ denote the set of equivalence classes
of n-dimensional complex representations
of the Weil-Deligne group $WD(F)$ on which Frobenius acts semisimply,
and let $\Cal{G}_0(n,F)$ denote its subset of irreducible representations.
In [H1] a map was constructed, $\pi \mapsto \sigma(\pi)$ from $\Cal{A}_0(n,F)$
to $\Cal{G}(n,F)$.  (The map was incorrectly normalized, however;
see the Erratum at the end of this introduction.)
This map enjoys a number of natural properties, some
of which are recalled below.  Using these properties, Henniart showed 
that this map
is a bijection with $\Cal{G}_0(n,F)$, and that this
bijection preserves Artin conductors ([He4], see [BHK]).
Bijections with most of these properties had already been constructed 
by Henniart [He2].
These properties do not suffice to determine the bijection uniquely, 
but Henniart proved
there is at most one such bijection compatible with the local epsilon factors
of pairs [He3].

More precisely, Langlands and Deligne have defined local epsilon factors
[L1,D2] for complex representations of the Weil group (more 
generally, of the Weil-Deligne
group); these local constants are compatible with the functional 
equations of complex
representations of the global Weil group.  The Langlands and Deligne 
constructions,
applied to the tensor product of two representations $\sigma_1 
\otimes \sigma_2$ of the Weil group,
yields the local epsilon factor of the pair $(\sigma_1, \sigma_2)$.
On the other hand, Jacquet, Piatetski-Shapiro,
and Shalika have defined local epsilon factors for pairs $(\pi_1, \pi_2)$, with
$\pi_1$ an admissible irreducible representation of $GL(n_i,F)$, $i = 
1, 2$; these
local constants are compatible with the functional equations of 
automorphic L-functions
of $GL(n_1,F) \times GL(n_2,F)$.  Henniart's
theorem is that, for any $n$, there is at most one
family of bijections $\Cal{A}_0(m,F) \leftrightarrow \Cal{G}_0(m,F)$, 
with $m \leq n$,
preserving local epsilon factors of pairs.  The local Langlands 
conjecture for $GL(n,F)$ can
thus be formulated as the assertion that a family of bijections 
preserving local
epsilon factors does indeed exist.  The main theorem of
the present article (Theorem 3.2) states that this is the case for $n < p$.

The correspondence
of [H1] is constructed globally.  One realizes the supercuspidal representation
$\pi$ as the local component at a place $v$ of a cohomological 
automorphic representation $\Pi$ of
$GL(n,E)$, where $E$ is a CM field of a certain type, with $E_v \cong 
F$.  In particular, one can
arrange that the global $L$-function of $\Pi$ is at almost all places the
$L$-function of a compatible family $\sigma(\Pi)$ of $\lambda$-adic 
representations
of $\GalE$, and one
takes $\sigma(\pi)$ to be the representation of a decomposition group 
at $v$.  Note
that it is not obvious that this correspondence is independent of $\ell$.
However, it is compatible with the Langlands correspondence
at unramified places, by construction.  It is therefore natural to expect that
$\pi \mapsto \sigma(\pi)$ ``is" the local Langlands correspondence.

Using a somewhat similar global approach, Laumon, Rapoport, and Stuhler proved
the local Langlands conjectures for local fields of positive characteristic.
To compare local epsilon factors, they followed an approach
introduced by Henniart [He1] and used the fact, observed by Deligne [D2],
that, in the case of function fields, the functional equation of the 
$L$-function
of a compatible family $\sigma(\Pi)$ of $\lambda$-adic representations
is consistent with the local epsilon factors of Langlands and Deligne.
This approach breaks down for number fields, since in general one 
doesn't know how
to prove the functional equation of the $L$-function of a compatible 
$\lambda$-adic
family, except by showing that the $L$-function is automorphic.  The 
exception, of
course, is when the compatible family is associated to a complex representation
of the Weil group.

In general,
an irreducible representation $\sigma$ of the local Weil group $W_F$ cannot be
globalized to a complex representation of the Weil group of a number field
whose associated compatible family of $\lambda$-adic representations
can be obtained from cohomological automorphic representations of $GL(n)$.
Specifically, as Henniart pointed out to me, a primitive 
representation (one not
induced from a proper subgroup)
cannot be so realized.  However, when $\sigma$ is monomial -- i.e., 
induced from
a character $\chi$ of a subgroup $H \subset W_F$ of finite index -- 
there is no obstacle in
principle to such a globalization.  The present paper constructs
globalizations of the desired sort when the fixed field $F'$ of $H$ 
is tamely ramified
over $F$.  In other words, we construct an extension $E'/E$ of global 
CM fields,
a place $v$ of $E$ with $E_v \cong F$, $E'_v \cong F'$, and an 
algebraic Hecke character
$\bc$ of $E'$ with local component isomorphic to $\chi$, (up to unramified
twist) such that $L(s,\bc)$ is equal
to the $L$-function of a cohomological cuspidal automorphic 
representation of $GL(n,E)$,
where $n = [F':F] = [E':E]$.  Here and below equality
of $L$-functions is understood to mean equality of Euler products, in this
case over the primes of $E$.

When $F'$ is a Galois extension of $F$ this is provided
by the global automorphic induction constructed by Arthur and Clozel in [AC]
(in the case of interest, the article [K] of Kazhdan suffices).
In this case, our main theorem is due to Henniart [He1].  The
difficulty is thus to construct automorphic induction
for non-Galois extensions, in enough cases to cover the local data of interest
(Lemma 1.6, Proposition 2.4).
Roughly speaking, we do this by working with a family of automorphic
representations to which, as Clozel has shown [C2], one can associate 
compatible families of
$\lambda$-adic Galois representations.  The Galois
representations are used to ``rigidify" the automorphic data.

It has been known for some time that every $\sigma \in 
\Cal{G}_0(n,F)$ is monomial
when $n$ and $p$ are relatively prime (cf. [KZ], [CH]; we use the reference
[KZ] for convenience).  In the literature, this is
referred to as the {\bf tame case}.  Using difficult calculations with Gauss
sums, Bushnell and Fr\"ohlich constructed bijections from $\Cal{G}_0(n,F)$
to $\Cal{A}_0(n,F)$ in the tame case, preserving $\epsilon$-factors 
for the standard
$L$-functions of $GL(n) \times GL(1)$ [BF].  [These bijections were based
on the internal structure of the multiplicative group of the division 
algebra, whose
; see also [M]).]  They also
observed that this property does not characterize such bijections uniquely.
By adapting the global methods of [He1]
to the case of non-Galois automorphic induction, we are able to show that the
bijections $\pi \mapsto \sigma(\pi)$ preserve $\epsilon$-factors for pairs of
irreducible Weil-group representations of degree prime to $p$.  When 
$n < p$, this
suffices, by [He3], to characterize the local correspondence uniquely.

It should be pointed out that the correspondence in [H1] is constructed
on $\ell$-adic cohomology, for $\ell \neq p$.  In principle it is 
possible that the
correspondence depends on $\ell$.  At least for $n < p$, the independence of
$\ell$ is a consequence of our main theorem.

More generally, we can realize any monomial representation of a Weil 
group, up to
unramified twist, as the local constituent of the compatible family of
$\lambda$-adic representations attached to a cuspidal automorphic 
representation obtained (by
non-Galois automorphic induction) from an appropriate Hecke 
character.  By Brauer's
theorem, the local epsilon factors of monomial representations 
suffice to determine the local
factors in general.  It is thus likely that the local Langlands 
conjecture can be
obtained for all $n$ by a generalization of the methods of the present paper.
In \S 4 we reduce the local Langlands conjecture in general to a partial
generalization to $GL(n)$ of Carayol's work [Ca] on the local Galois 
correspondence for Hilbert
modular forms at places of bad reduction.  Generalizing Carayol's results
will require a more complete understanding of the bad reduction of the
Shimura varieties used by Clozel to construct his compatible 
$\lambda$-adic families,
extending the results of [H1] and [H2] recalled in Theorem 1.7.

In the recent article [BHK] of Bushnell, Henniart, and Kutzko, it is
proved that any correspondence with the properties of the one 
constructed in [H1]
preserves {\it conductors} of pairs.  Their approach, valid in all degrees,
uses the relation between the Plancherel formula and the fine structure
theory of representations of $GL(n,F)$ -- the theory of types, due to Bushnell
and Kutzko.

The present paper makes extensive use of Henniart's ideas on the local
Langlands correspondence, especially those summarized in the letter [He4]
(incorporated in [BHK]).  Henniart's generous advice helped me to clarify
my ideas and spared me the consequences of a number of potentially 
embarrassing errors.
I also thank Clozel and Rohrlich for helpful comments.

\bigskip
\bigskip

\centerline{\bf Erratum to [H1]}
\medskip

The normalization of the correspondence $\pi \mapsto \sigma(\pi)$ in [H1]
was incorrect in two ways.  In the first place,
the conventions for the Shimura variety in [H1] were inconsistent with the
conventions of Rapoport and Zink for $p$-adic uniformization.  The $h_{\Phi}$
defined on p. 89 of [loc. cit.] is in fact the complex conjugate of the
Shimura datum $\bar{h}_{\Phi}$ determined by the conditions on the 
bottom of p. 302 of [RZ].

In the second place, the calculation of the twisting character $\nu$ on
pp. 100-101 was inconsistent with the definition of $h_{\Phi}$.
For the sake of completeness, and for future reference, here is the
correct formula for the twisting character for the Shimura datum 
determined by $h_{\Phi}$.
As in [H1], we let $Z$ be the connected center of $G\Cal{G}$.  Define
$r_{\mu}$ to be the algebraic
character $R_{\Cal{K}/\QQ} \Bbb{G}_{m,\Cal{K}} \ra Z$ which on 
$\QQ$-rational points is the map
$$\Cal{K}^{\times} \ra  Z(\QQ);~~ a \ra N_{\Cal{K}/\Cal{K}_0} a.$$
With this definition of $r_{\mu}$, the formula in [H1]:
$$\nu(G\pi) = \xi^{-1}\circ r_{\mu}\cdot 
|\bullet|_{\bold{A}}^{\frac{n-1}{2}} \tag E.1$$
makes Theorem 2 of [loc. cit.] correct.  Recall that Theorem 2 was a 
restatement
of results on the L-function of the Shimura variety, due to Kottwitz, Clozel,
and Taylor.  Similarly, the twisting character $\bar{\nu}$ for the 
Shimura datum $\bar{h}_{\Phi}$
is given by formula (E.1), with $r_{\mu}$ replaced by $\bar{r}_{\mu}$, where
$$\bar{r}_{\mu}(a)  =  r_{\mu}(\iota(a)); \tag E.2$$
here $\iota$ denotes complex conjugation.

The replacement of $\bar{h}_{\Phi}$ by $h_{\Phi}$
means we have inadvertently calculated the local Galois 
representation, not at the
prime $v$ where the Shimura variety admits $p$-adic uniformization, but rather
at its complex conjugate $\bar{v}$.  To continue, it is simplest (for 
reasons that
will become clear in the Appendix) to return to the conventions of 
Rapoport and Zink.  Define
$\bar{\nu}$ as above.  With the Shimura datum $(G,\overline{X}_{n-1})$,
where $\overline{X}_{n-1}$ is defined in terms of $\bar{h}_{\Phi}$, 
Theorem 2 of [loc. cit.]
becomes
$$\aligned
L^T(s,H^{n-1}(\overline{\Bbb{S}}_N,\Qbar_{\ell})[G\pi]) &=
L^T(s - \frac{n-1}{2},\Pi_{\Cal{K}}^*,St, \bar{\nu}_0(G\pi))^m \\
&= L^T(s,\Pi_{\Cal{K}}^*,St, \bar{\nu}(G\pi))^m.
\endaligned \tag E.3$$
Here $m$ is a multiplicity that plays no role in the present 
discussion.  The presence
of $\Pi_{\Cal{K}}^*$ rather than $\Pi_{\Cal{K}}$ in the formula differs from
the conventions of much of the literature but is consistent with the 
considerations
discussed on pp. 82-83 of [H3].

Reviewing the constructions in \S 3 of [H1], we see that
the local correspondence $\pi \mapsto \sigma(\pi)$ is
compatible with the global correspondence on cohomology (cf. Theorem 
1.7, below)
if we define  $\sigma(\pi_v)$ by
$$\sigma(\pi_v) = [\tilde{\sigma}(\pi_v)\otimes 
\bar{\nu}(G\pi)_v^{-1}]^*  \tag E.4$$
Here $*$ denotes contragredient and $\tilde{\sigma}(\pi_v)$ is defined
as in [H1]:
$$
\tilde{\sigma}(\pi_v) =
[Hom_{GJ}(H^{n-1}_c(\cbM_N,\Qbar_{\ell})_{SS(F)}, \pi_p) \otimes 
GJL(\pi_v)^{*}]^{GG}.  \tag E.5$$
With this definition, the global arguments in [H1] that show 
compatibility of the
correspondence $\pi \mapsto \sigma(\pi)$ with cyclic base change, 
automorphic induction,
local abelian class field theory, and so on, remain (or become) correct.

\bigskip
\bigskip

\centerline{\bf Notation}
\medskip

Let $G$ be a reductive group over the number field $E$.  For any
place $v$ of $E$, we write $G_v = G(E_v)$.  Let
$G_{\infty} = \prod G(E_{\sigma})$, the product taken over the set of
archimedean places of $G$; thus $E_{\infty}$ denotes the product
of the archimedean completions of $E$.  We also define $E_{\ad}$ 
(resp. $E_{\ad^f}$) to be the
adeles (resp. finite adeles) of $E$; if $S$ is a finite
set of finite primes of $E$ we let $E_{\ad^{f,S}}$ denote the ring of 
finite adeles
with entry $0$ at all places of $S$.  We let
$\Cal{A}(G)$ denote the space of automorphic forms on $G(E)\backslash 
G(E_{\ad})$,
(briefly:  automorphic forms on $G$), relative to an implicit choice of
maximal compact subgroup $K \subset G_{\infty}.$
We let $\Cal{A}_0(G) \subset \Cal{A}(G)$ be the space of cusp forms.
Let $Z_G$ denote the center of $G$, and let $\xi$ be a Hecke character
of $Z_G(E)\backslash Z_G(E_{\ad})$.
We let $\Cal{A}_0(G,\xi) \subset \Cal{A}_0(G)$ denote the subspace
of forms $f$ such that
$$f(zg) = \xi(z)f(g),~ z \in Z_G(E_{\ad}),~ g \in G(E_{\ad}).$$

We let $\Cal{T}(G)$ denote the set of automorphic representations
of $G$, which we take in the sense of irreducible admissible
representations of $G(\ad)$ that occur as subspaces of $\Cal{A}(G)$.
We define $\Cal{T}_0(G)$ to be the set of
cuspidal automorphic representations, and
$\Cal{T}_0(G,\xi)$ the set of cuspidal automorphic representations
with central character $\xi$.  If $S$ is a finite set of places
of $E$, let $\Cal{T}(G,S)$ denote the set of
automorphic representations of $G$ unramified outside $S$,
and define $\Cal{T}_0(G,S)$ and $\Cal{T}_0(G,S,\xi)$ in the obvious way.

If $\pi \in \Cal{T}(G)$, we let $\check{\pi}$ denote its contragredient.

Most frequently, we regard $G$ as a group over $\QQ$ by restriction
of scalars, and write $G(\QQ)$ and $G(\ad)$ instead of $G(E)$ and $G(E_{\ad})$.

By a Hecke character of the number field $E$ we mean a continuous 
complex character
of the idele class group $E_{\ad}^{\times}/E^{\times}$.  If $E'/E$ is 
a finite extension
of local or global fields, we let $N_{E'/E}$ and $Tr_{E'/E}$ denote the norm
and trace maps, respectively, from $E'$ to $E$.

Let $\Gamma$ be a group and let $H$ be a subgroup
of $\Gamma$ of finite index.  If $\sigma$ is a finite-dimensional
representation of $H$ over the coefficient field $k$, we let
$Ind_H^{\Gamma} ~\sigma$ denote the representation of $\Gamma$ 
induced from $\sigma$.
If $\Cal{L}/L$ is a (possibly infinite) Galois extension of fields,
$\Gamma$ is the Galois group $Gal(\Cal{L}/L)$, and $L' \subset 
\Cal{L}$ is the fixed
field of the (open and closed) subgroup $H$, we also write 
$Ind_{L'/L} \sigma$ for
$Ind_H^{\Gamma} ~\sigma$.
When $L$ is a local field and $\Gamma = W_L$ is its Weil group, $H = W_{L'}$
the Weil group of $L'$, then we again write $Ind_{L'/L} ~\sigma$ for
the induced representation of $\Gamma$.
Let $K \subset G$ be a second subgroup of finite index,
and let $A \subset G$ be a set of representatives of the double
cosets $H\backslash G/K$.  The {\bf Mackey constituents} of $Ind_H^G ~\sigma$,
restricted to $K$, are then the representations $Ind_{aHa^{-1}\cap 
K}^K ~a(\sigma)$,
for $a \in A$.  Here $a(\sigma)$ denotes the representation of 
$aHa^{-1}$ obtained
from $\sigma$ via the canonical isomorphism $aHa^{-1} \isoarrow H$.  Up to
isomorphism, the Mackey constituents do not depend on the set $A$ of double
coset representatives, and Mackey's Theorem is the isomorphism
$$Ind_H^G~ \sigma \isoarrow \oplus_{a \in A} Ind_{aHa^{-1}\cap K}^K 
~a(\sigma).$$

\bigskip
\bigskip
\centerline{\bf 1. Representations of $GL(n)$ and $\ell$-adic representations.}

\medskip

In the present section we let $E$ be an arbitrary number field and $G 
= GL(n)_E$,
viewed alternatively as a reductive group over $E$ or over $\QQ$, by 
restriction
of scalars.  For $\pi \in \CT(G)$ we let $L(s,\pi)$ denote the standard
(principal) $L$-function of $\pi$, as in [GJ], with the archimedean (Gamma)
factors excluded.  If $v$ is a finite place
of $E$ we let $L_v(s,\pi)$ be the local Euler factor at $v$ of $L(s,\pi)$, and
for any finite set $S$ of places of $E$ we write $L^S(s,\pi)$ for the partial
Euler product $\prod_{v \notin S} L_v(s,\pi)$, the product being taken over
finite places of $E$ not in $S$.

Let $\sigma = \{(\sigma_{\lambda},W_{\lambda})\}$ be a
compatible family of $\lambda$-adic representations of $\GalE$, where
$\lambda$ runs through the set of finite places of some number
field $L$, possibly with a finite set of $\lambda$ excluded, and
where $W_{\lambda}$ is a finite-dimensional vector space over $L_{\lambda}$.
We let $\ell(\lambda)$ denote the residue characteristic of $\lambda$.
Here by ``compatible" we mean that there is a finite set $S$ of
finite places of $E$ such that $\sigma_{\lambda}$ is
unramified for all $v$ outside $S \cup \text{ primes of residue 
characteristic }
\ell(\lambda)$ and that the characteristic polynomials
$P_{v,\lambda}(T) = det(1 - \sigma_{\lambda}(Frob_v)T)$ of
geometric Frobenius $Frob_v$ at $v$ have coefficients in $L$;
it is assumed that $P_{v,\lambda}(T) = P_{v,\lambda'}(T)$ as polynomials
in $L(T)$ for distinct primes $\lambda$ and $\lambda'$ of $L$, of
residue characteristic different from that of $v$, when $v \notin S$.  In
practice we will assume $\sigma$ to be the sum of representations
$\sigma_i$, with each $\sigma_i$ pure of some fixed weight $w$.  Thus
for $v \notin S$ the eigenvalues of $\sigma_i(Frob_v)$ have complex absolute
values $(Nv)^{\frac{w}{2}}$, for any complex embedding
of the number field the eigenvalues generate, with $Nv$ the cardinality of
the residue field of $v$. Let
$L^S(s,\sigma) = \prod_{v \notin S} L_v(s,\sigma)$ be the partial $L$-function
attached to the compatible family
$\sigma$.  We give $L^S(s,\sigma)$ the unitary (Langlands) normalization:
$$L_v(s,\sigma) = \prod_{\beta}(1 - \frac{\beta}{|\beta|}Nv^{-s})^{-1},$$
where $\beta$ runs over the set of reciprocal roots of $P_{v,\lambda}$
(any $\lambda$).  By our purity hypothesis the absolute values $|\beta|$
are well defined.

It should be borne in mind that $L^S(s,\sigma)$ in the unitary
normalization is not itself the (partial) $L$-function attached to a Galois
representation unless $w$ is even.  However, the unitary normalization is
convenient when applying base change.  In any case, $L^S(s-\frac{w}{2},\sigma)$
is the $L$-function attached to $\sigma$ in the arithmetic normalization.

\proclaim{Definition 1.1} The compatible family $\sigma$ is {\it 
weakly associated}
to $\pi \in \CT(G)$ if $L^{S'}(s,\sigma) = L^{S'}(s,\pi)$, as Euler products
over places of $E$, for some finite
set $S'$ containing $S$.  A compatible family $\sigma$ of dimension 
$n$ is called
{\it automorphic} if it is weakly associated to some $\pi \in \CT(G)$.
We let $\C(n,E) \subset \CT(G)$ denote the set of
$\pi \in \CT(G)$ for which there exists a compatible family $\sigma$
as above, weakly associated to $\pi$.
\endproclaim

Let $F$ be a $p$-adic field.  There is a one-to-one correspondence
between spherical representations $\tau$ of $GL(n,F)$ and unramified 
$n$-dimensional
completely decomposable complex representations of the Weil group 
$W_F$.  We denote
this correspondence $\tau \mapsto \sigma(\tau)$.  The relation
$L^{S'}(s,\sigma) = L^{S'}(s,\pi)$ of the preceding definition is 
equivalent to the
assertion that $\sigma_v$ is equivalent to $\sigma(\pi_v)$ for all 
finite $v \notin S'$.
Here $\sigma_v$ is the restriction of the compatible system $\sigma$ 
to a decomposition
group of $v$ in $\GalE$; $\pi_v$ is the $v$-component of $\pi$, an irreducible
admissible representation of $GL(n,E_v)$; and the compatible system $\sigma_v$
of unramified $\lambda$-adic representations of $\GalEv$ is identified with
a complex representation of $W_{E_v}$ in the usual way.  Of course, the local
representations $\sigma(\pi_v)$ can be attached to $\pi$ at all 
unramified places.

To each $\pi \in \CT(G)$ we may associate its cuspidal spectrum 
$Cusp(\pi)$.  This is
an unordered set of pairs $(n_i,\pi_i)$, $i = i,\dots, r$,
where $n_i$ is a positive integer and
$\pi_i \in \CT_0(GL(n_i,E))$, such that $n = \sum_i n_i$ and $\pi$ is a
subquotient of the induced representation
$Ind_P^G(\pi_1 \otimes \pi_2 \otimes \cdots \otimes \pi_r)$; here
$P \subset G$ is the standard parabolic subgroup attached to the
partition $n = \sum_i n_i$.
If $\sigma$ and $\pi$ are weakly associated, as above, then $\sigma$ 
and $Cusp(\pi)$
determine each other uniquely, by Chebotarev density and the 
classification theorem
of Jacquet-Shalika [JS], respectively.
Let $Isob(n,E) \subset \CT(G)$ denote the set of isobaric automorphic 
representations
of $G$ [L2, C1].  If $\pi \in Isob(n,E)$ then it is uniquely determined
by $Cusp(\pi)$, hence by a weakly associated $\sigma$.

We let $Reg(n,E) \subset Isob(n,E)$ denote the set of isobaric representations
$\pi$ such that, for each $\pi_i \in Cusp(\pi)$, the archimedean component
$\pi_{i,\infty}$ is of cohomological type [C1, Def. 3.12];
equivalently, $\pi_{i,\infty}$ has regular infinitesimal character. 
The following theorem
is mostly due to Clozel [C2]:

\proclaim{Theorem 1.2}  Let $E$ be a number field and let $\pi \in Reg(n,E)$.
Let $(n_i,\pi_i)$ be its cuspidal spectrum.  Suppose

(i) For each $i$ there exists a finite place $v(i)$ of $E$ such that the
local constituent $\pi_{i,v(i)}$ is square-integrable; let $p_i \in \QQ$ be the
rational prime divisible by $v(i)$, for each $i$.

(ii)  $E$ is of the form $E_0\cdot K_0$, with $E_0$ totally real
and $K_0$ imaginary quadratic and, for each $i$,
$\pi_i \simeq \check{\pi}_i\circ c$, where
$c: E \isoarrow E$ is the non-trivial element of $Gal(E/E_0)$.

We also assume that each of the primes $p_i$ splits in $K_0$.
Then $\pi \in \C(n,E)$.
\endproclaim

We let $\CU(n,E) \subset \C(n,E)$ denote the set of $\pi$ verifying 
the conditions
of Theorem 1.2.

Obviously it suffices to prove the theorem when $\pi$ is cuspidal.
Under a slightly strengthened version of assumption (i) (depending on
$n~\pmod{4}$), and (ii), the existence of a compatible system
$\sigma$ weakly associated to cuspidal $\pi$ was essentially proved by Clozel,
using the determination by Kottwitz of the zeta-functions of the Shimura
varieties to be discussed below, and Clozel's work on stable base change.
However, Clozel's argument only yielded $\sigma$ such
that $L^{S'}(s,\sigma) = L^{S'}(s,\pi)^m$, for some undetermined
multiplicity $m$.  The existence of $\sigma$ weakly associated
to $\pi$ was then deduced by Taylor ([T], unpublished; see [H1,\S 3]).
The sufficiency of (i) in the cases of bad parity was observed by 
Blasius, using
an argument based on quadratic base change (unpublished, but see [Bl, 4.7];
an alternative argument can be found in [C3], Theorem 2.6).

In (ii) one can take $E$ to be an arbitrary CM field, but then one only
obtains Galois representations over a certain reflex field containing
$E$ as a proper subfield, in general; cf. [C2, Th\'eor\`eme 5.3].  For
$E$ as in (ii) one verifies immediately, as in [BR, p. 66],
that the reflex field is just $E$ itself.

We let $\CU(n,E) \subset \C(n,E)$ denote the set of $\pi$ verifying 
the conditions
of Theorem 1.2; it is only defined for $E$ as in (ii).

\proclaim{Definition 1.3}  Let $E'$ be a finite extension of $E$.  Let
$\pi \in Isob(n,E)$ and $\pi' \in Isob(m,E')$, for some positive integers
$n$ and $m$.  Let $S$ be a finite set of finite places of $E$ and let $S'$
be the set of places of $E'$ above $S$.  We assume $S$ (resp. $S'$) 
contains all places
at which $\pi$ (resp. $\pi'$) is ramified.

(a)  Suppose $n = m$.  We say $\pi'$ is a {\bf weak base change} of 
$\pi$, or that
$\pi$ is a {\bf weak descent} of $\pi'$, relative to $S$, if, for any 
finite place
$v$ of $E$, $v \notin S$, and for any place $v'$ of $E'$ dividing $v$, we have
$$\sigma(\pi'_{v'}) = \sigma(\pi_v)|_{\GalEvp}.$$

(b)  Suppose $n = m[E':E]$.  We say $\pi$ is a {\bf weak automorphic 
induction} of
$\pi'$, relative to $S$, if, for any finite place
$v$ of $E$, $v \notin S$, and for any place $v'$ of $E'$ dividing $v$, we have
$$\sigma(\pi_v) = Ind_{E'_v/E_v} \sigma(\pi'_{v'}).$$
\endproclaim

Here $Ind_{E'_v/E_v}$ denotes induction from $\GalEvp$ to $\GalEv$,
as in the notation section.  Condition (b) is equivalent
to the condition that $L^S(s,\pi) = L^{S'}(s,\pi')$ (as partial Euler
products).  We say that $\pi'$ is
a weak base change of $\pi$ if it is a weak base change relative to 
some $S$, and
similarly for descent and automorphic induction.

If $E'$ is a solvable extension of $E$, Arthur and Clozel have proved 
the existence
of the weak base change map $BC_{E'/E}: Isob(n,E) \rightarrow 
Isob(n,E')$ and the
weak automorphic induction map $AI_{E'/E}: Isob(m,E') \rightarrow 
Isob(m[E':E],E)$,
relative to the sets $S$ and $S'$ of ramified places [AC, Theorems 
III.4.2, III.5.1,
III.6.2].  In fact, Arthur and
Clozel construct canonical local base change maps and show that their global
base change is compatible with local base change at all places.
(For our applications, the results of Kazhdan on cyclic automorphic 
induction of
characters are largely sufficient [K].)

\noindent (1.4) If $E'/E$ is cyclic of prime degree $\ell$, then
Arthur and Clozel determine the image and fiber of the base change
map [AC,Chapter 3].  Let $\alpha$ be a generator of $Gal(E'/E)$.
We denote the action of $Gal(E'/E)$ on $Isob(n,E')$ by $\pi' \mapsto {}^g\pi'$.
Then $\pi' \in Isob(n,E')$ is in the image of $BC_{E'/E}$ if and only if
$\pi' \cong {}^{\alpha}\pi'$.  The fibers of $BC_{E'/E}$
are completely determined by the following rules:

(a) If $\pi'$ is cuspidal and
$\pi' \simeq {}^{\alpha}\pi'$ then $\pi'$ is the base change of
precisely $[E':E]$ representations $\pi \in \CT_0(GL(n)_E)$,
all twists of one of them by powers of the class field character 
associated to $E'/E$.

(b) Suppose $m = {\ell}^{-1}\cdot n$ and suppose $Cusp(\pi') =
\{(m,\pi_1),(m,\pi_2),\dots,(m,\pi_{\ell})\}$, with $\pi_{i+1} \cong 
{}^{\alpha}\pi_i$
for $i = 1, \dots, \ell - 1$.  If $\pi_1 \ncong {}^{\alpha}\pi_1$ then $\pi'$
is the base change of exactly one $\pi \in Isob(n,E)$; moreover, 
$\pi$ is cuspidal.

The following lemma is an easy consequence of this description.

\proclaim{Lemma 1.5}  Let $\pi' \in Isob(n,E')$, and suppose $\pi' = 
BC_{E'/E}(\pi)$,
for some $\pi \in Isob(n,E)$.  Suppose $\pi \in \C(n,E)$.  Then $\pi' 
\in \C(n,E')$.
Moreover, any other weak descent $\pi_0 \in Isob(n,E)$ of $\pi'$ is also
in $\C(n,E)$.
\endproclaim

\demo{Proof}  If $\sigma$ is weakly associated to $\pi$, then $\sigma_{\GalEp}$
is weakly associated to $\pi'$.  The second assertion easily reduces 
to the case
in which $\pi$ is cuspidal, and then it follows from (a) and (b) above, since
$\C(n,E)$ is stable with respect to twists by characters of finite order.
\enddemo

Henniart and Herb have constructed canonical local automorphic induction maps
(for cyclic extensions of prime degree, hence for solvable extensions)
and shown that they are compatible at all places with the map 
$AI_{E'/E}$ of Arthur
and Clozel [HH].  The local automorphic induction is determined by certain
character relations.  We are interested in extending
the map $AI_{E'/E}$ to the case of an extension $E'/E$ which is not Galois but
whose Galois closure $\tilde{E}'$ is solvable over $E$.  Assuming 
$\pi' \in \CU(1,E')$ -- in particular,
$\pi'$ is an algebraic Hecke character -- then under additional regularity
hypotheses relative to the extension $E'/E$ one can construct a representation
$AI_{E'/E}(\pi') \in \CU([E'/E],E)$ whose weakly associated compatible family
is obtained by usual induction from the compatible family of one-dimensional
representations associated to $\pi'$.  In other words, $AI_{E'/E}(\pi')$ is a
weak automorphic induction of $\pi'$.

This is a vague formulation of a
general principle that will be discussed further in \S 4.  At present,
we restrict our attention to the simplest case, in which
$Gal(\tilde{E}'/E)$ is a semi-direct product of cyclic groups:
$$ \Gamma = Gal(\tilde{E}'/E) \cong A \rtimes T,$$
with $E'$ the fixed field of the non-normal subgroup $T$.  Moreover,
we assume we are in the situation of (ii) of Theorem 1.2:  there is 
an extension
$E'_0/E_0$ of totally real fields, with Galois closure $\tilde{E}'_0$, such
that $Gal(\tilde{E}'_0/E_0) \cong A \rtimes T$, with $E'_0$ the fixed field of
$T$, and an imaginary quadratic field $K_0$, such that $E = E_0\cdot K_0$,
$E' = E'_0 \cdot K_0$, etc.  We let $E^u \subset \tilde{E}'$ denote 
the fixed field
of $A$; it is a cyclic extension of $E$ with Galois group $T$.  We 
let $c$ denote
complex conjugation on any subfield of $\tilde{E}'$, and let $n = 
[E':E] = |A|$.

\proclaim{Lemma 1.6}  Let $\chi$ be an algebraic Hecke character of 
$E'$ such that
$\chi = \chi^{-1} \circ c$.  Let $\tilde{\chi} = \chi\circ 
N_{\tilde{E}'/E'}$ denote
the base change of $\chi$ to $\tilde{E}'$.  Suppose
$\pi^u = AI_{\tilde{E}'/E^u} \tilde{\chi} \in \CU(n,E^u)$ and is 
cuspidal, and let $v^u$ be a finite place
of $E^u$ at which $\pi^u$ is square-integrable; the existence of
such a place is provided by condition (i) of Theorem 1.2.  Suppose 
$v^u$ divides
the place $v$ of $E$ and is the only divisor of $v$.  Finally, suppose
the compatible family of $\lambda$-adic representations of $Gal(\Qbar/E^u)$
weakly associated to $\pi^u$ is irreducible.  Then there is a cuspidal
automorphic representation $\pi \in \CU(n,E)$ which is a weak automorphic
induction of $\chi$.
\endproclaim

\demo{Proof}  We let $\sigma'$ denote the compatible family of $\lambda$-adic
representations associated to $\chi$, and let $\tilde{\sigma}$ and 
$\sigma^u$ be
the compatible families (weakly) associated to $\tilde{\chi}$ and
$\pi^u$, respectively.  Then $\tilde{\sigma}$ is the restriction of 
$\sigma'$ to
$Gal(\Qbar/\tilde{E}')$ and $\sigma^u$ is the induction of $\tilde{\sigma}$ to
$Gal(\Qbar/E^u)$.  It follows that $\sigma^u$ is isomorphic to its conjugates
with respect to $T = Gal(E^u/E)$.  Thus $\pi^u$ is isomorphic at almost all
places to ${}^t(\pi^u)$ for any $t \in T$.  By strong multiplicity one,
$\pi^u$ is thus isomorphic to its $T$-conjugates, hence descends to a
cuspidal automorphic representation $\pi$ of $GL(n,E)$.  More 
precisely, as in (1.4.a) and
letting $m = |T|$, there is a set $\Cal{P}$ of $m$ distinct cuspidal automorphic
representations $\pi$ of $GL(n,E)$ that base change to $\pi^u$.  Let 
$\xi$ be a faithful
character of $T$.  Fixing $\pi \in \Cal{P}$,
the elements of $\Cal{P}$ are all of the form $\pi \otimes \xi^i \circ \det$,
$i = 1, \dots, m$, where $\xi$ is identified via class field theory
with a finite Hecke character of $E^{\times}$ trivial on the norms from $E^u$.

In order to complete the proof, it suffices to show that one, hence every
$\pi \in \Cal{P}$, belongs to $\C(n,E)$.  Indeed, assume $\pi \in \C(n,E)$, and
let $\sigma$ denote the weakly associated compatible family.  By the properties
of base change, $\sigma$ is an
extension of $\sigma^u$ to $Gal(\Qbar/E)$.  By hypothesis $\sigma^u$ 
is irreducible,
hence there are precisely $m$ such extensions, each obtained by 
twisting $\sigma$
by a power of $\xi$.  Thus every such extension is weakly associated to
exactly one element of $\Cal{P}$.  On the other hand, let $\sigma_1$ 
denote the compatible
family of $\lambda$-adic representations of $Gal(\Qbar/E)$ obtained by inducing
$\sigma'$.  It is elementary to verify that $\sigma_1$ is an 
extension of $\sigma^u$.
Without loss of generality, we may thus assume $\sigma_1 = \sigma$, 
and $\pi$ is
then a weak automorphic induction of $\chi$.  Since $\chi$ satisfies
$\chi = \chi^{-1} \circ c$, it follows that $\sigma' \circ c$ is isomorphic
to the contragredient of $\sigma'$.  Thus $\sigma$ has the same property.
Since $\sigma$ is weakly associated to $\pi$,
it follows from strong multiplicity one that $\pi \simeq \check{\pi}\circ c$.
Thus $\pi$ actually belongs to $\CU(n,E)$.

We first verify that $\pi \in Reg(n,E)$ and satisfies condition (i) 
of Theorem 1.2.
First, $\pi \in Reg(n,E)$ because $\pi^u \in Reg(n,E^u)$,
since for any complex place $w$ of $E$ we have $\pi_w \cong \pi_{w^u}$ for any
place $w^u$ of $E^u$ dividing $w$.  Next, condition (i) at the place $v$ is
automatic since $v^u$ is the only divisor of $v$ in $E^u$ and by the hypothesis
on $\pi^u$ at $v^u$.

Finally, we verify a weak analogue of (ii) of Theorem 1.2 that will 
be sufficient
for our purposes.  We recall that base change
from $E$ to $E^u$ commutes with complex conjugation and with the 
contragredient.
Thus $BC_{E^u/E}(\check{\pi}\circ c) = \check{\pi}^u\circ c = \pi^u$.
It then follows that $\check{\pi} \circ c \cong \pi \otimes \eta \circ det$,
where $\eta$ is a power of $\xi$.   The abelian extension $E^u/E$ is a lift
of a totally real abelian extension of $E_0$; thus $\eta = \eta_0 
\circ N_{E/E_0}$, for
some finite Hecke character $\eta_0$ of the ideles of $E_0$, trivial 
at the archimedean
places of $E_0$.   Let $\alpha$ be a finite Hecke character of 
$E^{\times}_{\ad}/E^{\times}$
that restricts to $\eta_0$ on the ideles of $E_0$.  Such an $\alpha$ exists:
indeed, $\eta_0$ is trivial at infinity, hence it suffices to 
construct a continuous
character of the compact group $E^{\times}_{\ad}/E^{\times}\cdot 
E_{\infty}^{\times}$
that restricts to $\eta_0$ on the closed subgroup
$E^{\times}_{0,\ad}/E^{\times}_0\cdot E_{0,\infty}^{\times}$, and this
is certainly possible.
Let $\pi_1 = \pi \otimes \alpha$,
where for brevity we write $\otimes \alpha$ instead of $\otimes 
\alpha \circ det$.
Then $\pi_1$ is still regular and still satisfies (i); moreover, we have
$$
\align
\check{\pi}_1 \circ c &= \check{\pi}\circ c \otimes \alpha^{-1} \circ c \\
&= \pi \otimes \eta_0 \circ N_{E/E_0} \cdot \alpha^{-1} \circ c \\
& = \pi_1;
\endalign$$
the last equality follows from the identity
$$\alpha(a) \cdot \alpha(a^c) = \alpha(N_{E/E_0}(a))
= \eta_0(N_{E/E_0}(a)).$$
Thus $\pi_1 \in \CU(n,E) \subset \C(n,E)$.  But $\C(n,E)$ is obviously
invariant under twists by finite Hecke characters, and we are done.

\enddemo

\noindent{\bf Remark.}  Under appropriate hypotheses, one can use 
similar arguments
to construct base change for the non-Galois extension $E'/E$, or more generally
in the setting of Proposition 4.8, below.

Let $F$ be a finite extension of $\QQ_p$.  In [H1] a map $\pi \mapsto 
\sigma(\pi)$
is constructed from
the set of equivalence classes of supercuspidal representations of $GL(n,F)$
to the set of $n$-dimensional representations of the Weil group $W_F$ of $F$.
It is shown that this map is compatible with isomorphisms $F \isoarrow F'$,
commutes with twists by characters of $F^{\times} \cong W^{ab}_F$, and takes
base change and local automorphic induction to restriction and induction
of Weil group representations, respectively.  Using the results and techniques
of [He2], Henniart showed that these conditions implied that $\pi 
\mapsto \sigma(\pi)$
is a bijection with the set of {\it irreducible}
representations and preserves Artin conductors [He 4, BHK].  It also 
commutes with taking
contragredients, though this was not stated explicitly.

We will need the following property of this bijection.

\proclaim{Theorem 1.7}  Let $E$ be a number field and let $\pi \in \CU(n,E)$.
Let $v$ be a finite place of $E$ and suppose the local component $\pi_v$ at
$v$ is supercuspidal.  Let $\sigma$ be the compatible family of $\lambda$-adic
representations of $\GalE$ weakly associated to $\pi$, and let 
$\sigma_v$ denote
the restriction of $\sigma$ to a decomposition group $\GalEv$ of $v$.  Then
$\sigma_v$ is equivalent to $\sigma(\pi_v)$.
\endproclaim

As noted in the introduction,
$\sigma(\pi_v)$ is really a family $\{\sigma(\pi_v)_{\ell}\}$ of 
representations
with coefficients in
$\Qbar_{\ell}$, for all $\ell \neq p$.  The theorem should be 
understood to mean
that $\sigma_{v,\ell}$ is equivalent to $\sigma(\pi_v)_{\ell}$ for each $\ell$.

The complete proof of this theorem will be given elsewhere [H2].  However, when
$\pi_{\infty}$ has cohomology with coefficients in the trivial representation
(i.e., has the same infinitesimal character as the trivial representation),
then this theorem is essentially proved in [H1].  Indeed, $\sigma(\pi_v)$
is defined to make Theorem 1.7 true in this case.  The work in [H1] 
is to prove,
using $p$-adic uniformization, that $\sigma(\pi_v)$ thus defined depends only
on the supercuspidal representation $\pi_v$ and not on the automorphic
representation $\pi$.  The case of cohomology with trivial coefficients
suffices for the applications in \S 3.  The arguments in \S 4 require
more general coefficient systems.  The reader may therefore prefer to
consider the results of \S 4 conditional, pending appearance of [H2].

It should also be noted that the choice of level subgroup
in [H1] (especially in the appendix) forces the local constituent of
$\pi$ at every place $w \neq v$ dividing $p$
to be a twist of the Steinberg representation.  However, this hypothesis
was only made in order to simplify notation and is irrelevant
to the final result.  In particular, the $p$-adic uniformization of
[H1,(A.11)] is valid with minor modifications for general
level subgroups at the primes $w$ as above.  At the referee's request,
we explain how to remove this hypothesis in the appendix.

Finally, the local correspondence in [H1] is constructed on $\ell$-adic
cohomology and therefore associates $n$-dimensional $W_F$-modules 
over $\Qbar_{\ell}$
to supercuspidal representations of $GL(n,F)$ with coefficients in 
$\Qbar_{\ell}$.
To obtain a candidate for the local Langlands correspondence we need to choose
an isomorphism between $\Qbar_{\ell}$ and $\CC$.  On the other hand, 
in the present
article we work with local components of induced representations from
algebraic Hecke characters $\chi$.  These take their values a priori 
in $\CC$, but
the local components can all be defined over the subfield $L(\chi) \subset \CC$
generated by the values of $\chi$ on the finite ideles.  Then $L(\chi)$ is a
number field and it is easy to see that the representations 
constructed by automorphic
induction are also defined over $L(\chi)$, up to twisting by a
half-integer power of the absolute value character (which may introduce
a square root of $p$).  Of course the association of compatible
families of $\ell$-adic representations to complex-valued Hecke
$L$-functions goes back to Weil.  The relation between complex and 
$\ell$-adic epsilon
factors is worked out in [D2,\S 6].

In particular, the number field $L$ that appears above Definition 1.1 
will be of
the form $L(\chi)$, and in particular is given as a subfield of 
$\CC$.  We could just
as well have adopted this point of view in [H1].  Indeed, the global 
cohomological automorphic
representation $G\pi$ used on p. 95 ff of [H1] to construct the local 
correspondence
can be realized over a number field $L(\pi)$ which is also a subfield 
of $\CC$.  By
varying $G\pi$ while keeping the local component $\pi_v$ fixed it 
should be possible
to define $\sigma(\pi_v)$ directly over a fixed finite extension of
the (number) field of definition of $\pi_v$ (with the square root of 
$p$ adjoined,
as above).  This finite extension should be given
by the field of definition of the representation of the 
representation $JL(\pi_v)$ of
the multiplicative group of the division algebra over $F$ with 
invariant $\frac{1}{n}$,
where the notation $JL(\bullet)$ denotes the Jacquet-Langlands 
correspondence, as in [H1].
In this way it may be possible to verify that the correspondence of 
[H1] is in fact independent
of $\ell$; however, I have not carried this out in detail.  Note that 
this argument would
depend on multiplicity one theorems for automorphic representations of unitary
similitude groups that are not presently in the literature.

\bigskip

\centerline{\bf 2. Tame representations of local Galois groups.}

\medskip

For any finite extension $\Phi$ of $\QQ_p$ we let $\G_{\Phi}$ denote 
$Gal(\Qbar_p/\Phi)$;
$G_{\Phi}^0$ denotes the inertia subgroup.
By a {\bf representation} of $\G_{\Phi}$ we will always mean a 
finite-dimensional
representation over an algebraically closed field of characteristic zero
that factors through a finite quotient.

Fix a finite extension $F$ of $\QQ_p$.
A representation of $\G_F$ is {\it tame} if
its restriction to the wild ramification subgroup $G_F^1$ decomposes 
as the sum of
one-dimensional representations; in other words, the restriction 
factors through
the abelianization of $G_F^1$.  If $(n,p) = 1$ then it is easy to see
that any irreducible $n$-dimensional representation of $G_F$ is tame 
(see, e.g. [KZ, 1.8]).
Only slightly more difficult to show is the fact that any irreducible tame
representation of $\G_F$ is {\bf monomial}, i.e. induced from a 
character of a subgroup
of the form $\G_{F'}$, for $F'$ a finite tamely ramified extension of 
$F$ [KZ, p. 104]
(cf. also [CH] and [D1], Proposition 3.1.4).

The characters of $\G_{F'}$ can be identified via class field theory 
with characters $\chi$ of
$(F')^{\times}$ of finite order; we will use the same notation to designate
the two associated characters.  Koch and Zink determine the pairs 
$(F'/F,\chi)$, with
$F'/F$ tamely ramified and $\chi$ a character of $\G_{F'}$,
such that $Ind_{F'/F} \chi$ is irreducible, and construct a bijection 
between the
set of classes of such pairs, under a natural equivalence relation, and the set
of equivalence classes of tame representations of $\G_F$ [KZ, Theorem 
1.8 and 3.1].
Here the notation $Ind_{F'/F}$ is as in the notation section.  For 
our purposes, it suffices
to note the following lemma:

\proclaim{Lemma 2.1} (Koch, Zink).  Suppose $F'/F$ is a tamely and 
totally ramified
finite extension and $\chi$ is a character of $\G_{F'}$ such that 
$Ind_{F'/F} \chi$ is
irreducible.  Let $\chi^1$ denote the restriction of $\chi$ to the 
wild ramification
subgroup $\G_F^1$.  Then $\G_{F'}$ is the stabilizer of $\chi^1$ in $\G_F$.
\endproclaim

A complete proof is given on p. 104 of [KZ], though the lemma is not 
stated as such.
Specifically, Koch and Zink begin by letting $\chi^1$ be any 
character of $\G_F^1$
contained in the restriction of $Ind_{F'/F} \chi$ to $\G_F^1$; 
obviously our chosen
$\chi^1$ fits that description.  They then define $K_1$ to be the 
fixed field of the
stabilizer of $\chi^1$.  Thus $F' \supset K_1 \supset F$.
Finally, they deduce that $F'$ is unramified over $K_1$, hence equals 
$K_1$ under
our hypotheses.

\proclaim{Corollary 2.2}  Under the hypotheses of Lemma 2.1, let 
$\tilde{F}$ denote
the Galois closure of $F'$ over $F$ and let $F^u \subset \tilde{F}$ denote the
maximal subextension unramified over $F$.  Let $\tilde{\chi}$ denote
the restriction of $\chi$ to $\G_{\tilde{F}}$.  Then the stabilizer 
of $\tilde{\chi}$ in
$Gal(\tilde{F}/F^u)$ is trivial, and $Ind_{\tilde{F}/F^u} 
\tilde{\chi}$ is irreducible.
\endproclaim

\demo{Proof}  The extension $F'/F$ is tame, hence $\G_{\tilde{F}}$ 
contains $\G_F^1$.
It thus follows from Lemma 2.1 that the stabilizer of $\tilde{\chi}$ 
in $Gal(\tilde{F}/F^u)$
is trivial.  The final assertion follows immediately, since 
$Gal(\tilde{F}/F^u)$ is cyclic.
\enddemo

We now return to the global setting of \S 1.  Let $E'/E$ be an 
extension of number fields
of degree $n$, as in the discussion preceding Lemma 1.6; in 
particular, it is assumed
that $E = E_0 \cdot K_0$ and $E' = E'_0 \cdot K_0$.  Define the 
fields $\tilde{E}$ and $E^u$,
so that $ \Gamma = Gal(\tilde{E}/E) \cong A \rtimes T$, with $E^u$ 
the fixed field of $A$
and $E'$ the fixed field of $T$.   We assume there is a finite place 
$v$ of $E_0$ such that $E_{0,v}$
is the $p$-adic field $F$ discussed above, and such that $v$ is 
divisible by exactly
one place $\tilde{v}$ of $\tilde{E}_0$.  Let $v'$ and $v^u$ be the 
corresponding
places of $E'_0$ and $E^u_0$, respectively.  We assume the rational 
prime $p$ splits in
$K_0$, so that $v$ splits in $E$ as the product of two primes $v_1$ 
and $v_2$; we define
the primes $v'_1$ and $v'_2$ of $E'$, $v^u_1$ and $v^u_2$ of $E^u$, 
and $\tilde{v}_1$
and $\tilde{v}_2$ of $\tilde{E}$ analogously, and identify $F$ with $E_{v_1}$.
We let $\tilde{F}$ denote the completion of $\tilde{E}$
at $\tilde{v}_1$ and define $F'$ and $F^u$ analogously; thus the decomposition
group $Gal(\tilde{F}/F)$ is isomorphic to $\Gamma$.   We assume 
$\tilde{F}$ is tamely
ramified over $F$.  Under our hypotheses, $\tilde{F}$ is necessarily 
unramified over $F'$.

\proclaim{Definition 2.3}  An algebraic Hecke character $\chi$ of 
$E'$ is called
$E'/E$-{\it regular} if, for any archimedean place $\tau$ of $E$, the 
components
of $\chi$ at the distinct places of $E'$ dividing $\tau$ are distinct.
\endproclaim

\proclaim{Proposition 2.4}  Let $\chi_v$ be a character of 
$(F')^{\times}  \cong W_{F'}$
with the property that $Ind_{F'/F} \chi_v$ is an irreducible representation of
$W_F$.  Let $\chi$ be an algebraic Hecke character of $E'$ with 
component $\chi_v$ at $v'_1$.
Suppose $\chi = \chi^{-1} \circ c$.  Suppose $\chi$ is 
$E'/E$-regular.  Then there is a cuspidal
automorphic representation $\pi \in \CU(n,E)$ which is a weak 
automorphic induction of $\chi$.
Moreover, the local component $\pi_{v_1}$ at $v_1$ is supercuspidal, 
and $\sigma(\pi_{v_1})$
is equivalent to $Ind_{F'/F} \chi_v$, where $\sigma(\bullet)$ is the 
correspondence of [H1],
discussed in \S 1.  Finally, at every archimedean place $\tau$ of $E$ we have
$$\pi_{\tau} \simeq \chi_{\tau_1} \boxplus \cdots \boxplus 
\chi_{\tau_n}, \tag 2.4.1$$
where $\boxplus$ denotes Langlands sum and $\{\tau_1,\dots, \tau_n\}$ 
is the set of
places of $E'$ above $\tau$.
\endproclaim

\demo{Proof}  We first assume $F'/F$ is totally ramified, and verify 
the conditions
of Lemma 1.6 in this case.  Define $\pi^u = AI_{\tilde{E}/E^u} 
\tilde{\chi}$ as in the
statement of Lemma 1.6.  I claim that the local component of $\pi^u$ 
at $v^u_1$ is
supercuspidal.  Indeed, there is an unramified character $\beta$ of 
$(F')^{\times}$
such that $\chi^*_v = \chi_v \otimes \beta$ is a character of finite 
order. Since
$$Ind_{F'/F}(\chi_v \otimes \beta_v) \cong Ind_{F'/F}(\chi_v) \otimes \beta_v$$
we see that
$Ind_{F'/F} \chi^*_v$ is irreducible.
By Corollary 2.2 the stabilizer of $\tilde{\chi}_v$ in the cyclic
group $Gal(\tilde{F}/F^u) \cong A$, which is also
the stabilizer of $\tilde{\chi^*}_v$, is trivial; thus 
$AI_{\tilde{F}/F^u} \tilde{\chi}_v$
is supercuspidal [HH, Proposition 5.5].  It follows that $\pi^u$ is 
cuspidal.  Moreover,
it follows as in the proof of Proposition 5 of [H1] that $\pi^u \cong 
\check{\pi}^u \circ c$.
Finally, we verify that $\pi^u \in Reg(n,E^u)$.  Note that the hypothesis of
$E'/E$ regularity implies that $\tilde{\chi}$ is 
$\tilde{E}/E^u$-regular, in the obvious
sense.  Let $\tau$ be an archimedean place of $E^u$ and let $\nu_j$, 
$j = 1, \dots, n$,
$n = [\tilde{E}/E^u] = |A|$,
denote the places of $\tilde{E}$ dividing $\tau$.  Then the local component
of $\pi^u$ at $\tau$ is
$$\tilde{\chi}_{\nu_1} \boxplus \tilde{\chi}_{\nu_2} \boxplus \cdots 
\boxplus \tilde{\chi}_{\nu_n}  \tag 2.4.2 $$
where $\boxplus$ denotes Langlands sum (cf [H1,(22)]).  Thus the 
regularity of $\pi^u$
is a consequence of the $\tilde{E}/E^u$-regularity of $\tilde{\chi}$.

It follows that $\pi^u \in \CU(n,E^u)$, and that its component at $v^u_1$,
which we denote $\pi^u_v$, is supercuspidal.
Let $\sigma^u$ denote the compatible family of $\lambda$-adic 
representations of $Gal(\Qbar/E^u)$
weakly associated to $\pi^u$.  By Theorem 1.7 the restriction
of $\sigma^u$ to a decomposition group of $v^u_1$ is equivalent to 
$\sigma(\pi_v)$.  By Henniart
[He4] this is an irreducible representation, hence $\sigma^u$ is {\it 
a fortiori} irreducible.
Of course $\sigma^u$ is equivalent to $Ind_{\tilde{E}/E^u} \tilde{\chi}$.

We have thus verified all the conditions of Lemma 1.6, from which the 
existence of
the weak automorphic induction $\pi$ of $\chi$ follows.  Moreover, by 
the construction
in Lemma 1.6, the base change to $F^u$ of the local component 
$\pi_{v_1}$ is the supercuspidal
representation $\pi^u_v$; hence $\pi_{v_1}$ is itself supercuspidal 
[AC, I.6; this fact
is proved but not stated during the proof of Lemma I.6.10].  Next,
letting $\sigma$ denote the compatible family weakly associated to 
$\pi$, it follows
from Definition 1.3 (b) that $\sigma(\pi_w) = Ind_{E'_w/E_w} 
\sigma(\pi'_{w'})$ for
almost all unramified places $w$, where the notation is as in Definition 1.3.
The assertion regarding $\sigma(\pi_{v_1})$ thus follows from 
Chebotarev density, as in the proof
of Proposition 5 of [H1].  Finally, (2.4.1) is an immediate 
consequence of (2.4.2).

This completes the proof of the proposition, provided $F'/F$ is 
totally ramified.
In general, let $F_1/F$ be the
maximal unramified extension contained in $F'$, and let $E_1$ be the 
corresponding
extension of $E$ contained in $E'$.  The irreducibility of
$Ind_{F'/F} \chi_v$ obviously implies irreducibility of $Ind_{F'/F_1} 
\chi_v$, and
$E'/E_1$-regularity of $\chi$ follows from $E'/E$-regularity.
Thus the above argument shows the existence of the weak automorphic induction
$\pi_1 \in \CU(n_1,E_1)$, with $n_1 = [F_1:F] = [E_1:F]$.  But $F_1/F$ is
unramified, hence cyclic.  It follows that $E_1$ is a cyclic extension of
$E$, and thus the automorphic induction $\pi = AI_{E_1/E}(\pi_1)$ is defined.
Obviously, $\pi$ is a weak automorphic induction of $\chi$, and it remains
to show that $\pi_{v_1}$ is supercuspidal, the final assertion following from
Chebotarev density as before.  Let
$v_{11}$ denote the prime of $E_1$ dividing $v_1$ and let $\pi_{11}$ denote
the local component of $\pi_1$ at $v_{11}$.  We know from the totally
ramified case that $\pi_{11}$ is supercuspidal, and that
$\sigma(\pi_{11})$ is equivalent to $Ind_{F'/F_1} \chi_v$.  It remains to show
that $\pi_{11}$ has trivial stabilizer in the cyclic group $Gal(F_1/F)$.  But
the local correspondence $\bullet \mapsto \sigma(\bullet)$ is 
$Gal(F_1/F)$-equivariant,
so it suffices to show that $Ind_{F'/F_1} \chi_v$ has trivial stabilizer
in $Gal(F_1/F)$.  This follows immediately from the irreducibility of
$Ind_{F'/F} \chi_v$.

\enddemo

\noindent{\bf Remark.}  The proof of Proposition 2.4 makes use
of Theorem 1.7, whose proof has not been published in general.  For the purpose
of comparing epsilon factors of pairs in the tame case,
it will suffice to consider $\chi$ for which the Langlands sum in 
(2.4.1) is a cohomological
representation with the infinitesimal character of the trivial representation,
as in [H1, (22)].  In that case, as we have already mentioned, 
Theorem 1.7 is proved
in [H1].

The representations constructed in Proposition 2.4 will be called 
{\it of monomial type}.
More general monomial type representations (with more general solvable 
Galois groups)
will be considered in \S 4.

\bigskip
\bigskip

\centerline{\bf 3.  Comparison of epsilon factors of pairs}
\medskip

The existence of a global extension $\tilde{E}/E$ adapted to $\tilde{F}/F$,
as in \S 2, is guaranteed by [D2, Lemma 4.13].  The fact that $E$ can be taken
in the form $E_0 \cdot K_0$ and so on follows from an easy approximation
argument.

For any finite place $w$ of $E'$ we let $U_w$ denote the unit 
subgroup of the multiplicative
group of the completion $E'_w$.
We let $\chi_v^0$ denote the restriction of $\chi_v$ to the unit 
subgroup $U_{F'}$
of $(F')^{\times}$.  Obviously, there exists a Hecke character $\gamma$ of $F'$
of finite order, unramified at $v'_2$ and equal to $\chi_v$ on 
$U_{v'_1} \cong U_{F'}$.
We let $\chi(1) = \gamma \cdot \gamma^{-1}\circ c$.
Then $\chi(1)$ restricts to $\chi_v$ on $U_{F'}$.

Next, for each archimedean (complex) place $\nu$ of $E'$ we fix a 
local character
$\chi_{\nu}$ of $\CC^{\times}$, trivial on $\RR^{\times}$, such that
$\chi_{\nu_1} \neq \chi_{\nu_2}$ whenever $\nu_1$ and $\nu_2$ restrict to the
same place on $E$.  Let
$$\chi_{\infty} = \prod_{\nu} \chi_{\nu}: 
(E'_{\infty})^{\times}/(E'_{0,\infty})^{\times}
\rightarrow \CC^{\times}$$
denote the corresponding infinity type.  Extend the character
$\chi_{\infty}$ trivially to a character $\delta$ of
$(E'_{\infty})^{\times} \times U_{v'_1} \times U_{v'_2})\cdot 
(E'_{0,\ad})^{\times}\cdot (E')^{\times}$,
viewed as a subgroup of the ideles of $E'$.  Let $\delta(1)$ denote 
any extension of $\delta$
to a Hecke character of $E'$, and let $\chi = \chi(1)\cdot \delta(1)$.  Then
$\chi$ is an $E'/E$-regular Hecke character of $E'$, whose restriction
to $(F')^{\times}$ equals $\chi_v$, up to an unramified twist. 
Moreover, $\chi$ is
trivial on $(E'_{0,\ad})^{\times}$, by construction, hence satisfies
$\chi = \chi^{-1} \circ c$.

Finally, we assume that every local component $\chi_{\nu}$ is of the form
$(z/\overline{z})^{\frac{a(\nu)}{2}}$, where the $a(\nu)$ are all 
integers of the
same parity.
Then if $a$ is any integer congruent to $a(\nu) \pmod{2}$
(any $\nu$), the product $\chi \cdot ||\bullet||^{\frac{a}{2}}$ is a 
Hecke character
of type $A_0$, in Weil's terminology (also called a motivic Hecke character).

In order to obtain cohomological representations
of $GL(n)$ by automorphic induction we will take $a(\nu) ~\equiv ~n-1 \pmod{2}$
We note that any $E'/E$-regular infinity type with the given parity 
condition can be used
in this construction.  If the infinity type is chosen as in 
[H1,(23)], where the present
$n$ replaces the $d$ of [H1] and $n$ of
[loc. cit.] is taken to equal $1$ -- the hypothesis of [loc. cit.] that
$E'/E$ is cyclic is irrelevant to the present construction -- we obtain
automorphic representations contributing to cohomology with trivial 
coefficients.

More generally, suppose $F_1/F$ and $F_2/F$ are two finite tame 
extensions; let $F'$ denote
their compositum and $\tilde{F'}$ its Galois closure.  Then we can 
find totally real global
fields
$$\tilde{E'}_0 \supset E'_0 = E_{1,0}\cdot E_{2,0} \supset E_0$$
with $Gal(\tilde{E'}/E) \cong Gal(\tilde{F}/F)$, as before.  We 
choose an imaginary quadratic
field $K_0$ in which $p$ splits and define $E = E_0\cdot K_0$, $E_1 = 
E_{1,0}\cdot K_0$,
and so on.  The above argument yields:

\proclaim{Lemma 3.1}  Let $F$ be a $p$-adic field and let $F_1$ and 
$F_2$ be finite
tame extensions of $F$ of degree $n_1$ and $n_2$, respectively.  Let 
$\chi_i$ be
characters of $G_{F_i}$,
such that $Ind_{F_i/F} \chi_i$ is irreducible, $i = 1, 2$.  Then there exists
a global CM field $E = E_0 \cdot K_0$, with $K_0$ imaginary quadratic and
$E_0$ totally real, and cuspidal automorphic representations $\Pi_i 
\in \CU(n_i,E)$, $i = 1, 2$,
with the following properties:

(a)  Both $\Pi_1$ and $\Pi_2$ are of monomial type;

(b)  There is a place $v_1$ of $E$ with $E_{v_1} \cong F$,
and the local constituents of $\Pi_1$ and $\Pi_2$ at $v_1$  are supercuspidal.

(c) Let $\Sigma_i$ be the compatible family of $\lambda$-adic representations
weakly associated to $\Pi_i$, $i = 1, 2$.  There are unramified characters
$\beta_i$ of $F_i^{\times}$, $i = 1, 2$, so that the restriction of 
$\Sigma_i$ to
a decomposition group of $v$ is equivalent to $(Ind_{F_i/F} \chi_i) 
\otimes \beta_i$,
$i = 1, 2$.
\endproclaim

We now come to the main theorem of this paper.

\proclaim{Theorem 3.2}  Let $n_1$ and $n_2$ be two positive integers 
prime to $p$.  Let
$\pi_1$ and $\pi_2$ be supercuspidal representations of $GL(n_1,F)$ 
and $GL(n_2,F)$,
respectively.  Let $\sigma_i = \sigma(\pi_i)$, $i = 1, 2$, be the corresponding
irreducible representations of $W_F$.  Then we have an equality of local
factors:
$$\epsilon(\pi_1 \times \pi_2,s,\psi,dx) = \epsilon(\sigma_1 \otimes 
\sigma_2, s, \psi, dx). \tag 3.2.1$$
Here $\psi$ is an arbitrary additive character of $F$, $dx$ is a Haar measure
on $F$, self-dual relative to $\psi$,
the local factor on the right is that defined in [JPS], and that on the
left is the one defined in [D2].

In particular, if $m \leq n < p$, and if $\pi_1$ and $\pi_2$ are irreducible
admissible representations of $GL(n,F)$ and $GL(m,F)$, respectively, then
the equality (3.2.1) holds, where $\sigma$ is extended to general irreducible
admissible representations by the procedure of [BZ].  The map
$\pi \mapsto \sigma(\pi)$ is the only bijection $\Cal{A}_0(n,F) 
\leftrightarrow \Cal{G}_0(n,F)$
with this property for $n < p$.
\endproclaim

\demo{Proof} By the results of Koch and Zink recalled in \S 2, we may 
assume $\sigma_i$
to be of the form $Ind_{F_i/F} \chi_i$, $i = 1, 2$ with $\chi_i$ as 
in Lemma 3.1.  Let
$\Pi_1$ and $\Pi_2$ be the cuspidal automorphic representations whose 
existence is
asserted in Lemma 3.1, and let $\Sigma_1$ and $\Sigma_2$ be the weakly
associated compatible $\lambda$-adic families.  Then we have an 
identity of partial $L$-functions
$$L^S(\Pi_1 \times \Pi_2, s) = L^S(\Sigma_1 \times \Sigma_2, s), \tag 3.2.2$$
where $S$ is a finite set of places.  Now Lemma 3.1 (a) implies that 
the right-hand
side of (3.2.2) is the partial $L$-function of a complex 
representation of the Weil group of
$E$.  Thus the completed L-function $L(\Sigma_1 \times \Sigma_2, s)$ 
satisfies a functional
equation of the form
$$L(\Sigma_1 \times \Sigma_2, s) =
\epsilon(\Sigma_1 \times \Sigma_2, s) L(\check{\Sigma}_1 \times 
\check{\Sigma}_2, 1-s), $$
where $\epsilon(\Sigma_1 \times \Sigma_2, s)$ is the product of the 
local factors of Langlands
and Deligne.  On the other hand, the left-hand side satisfies
$$L(\Pi_1 \times \Pi_2, s) =
\epsilon(\Pi_1 \times \Pi_2, s) L(\check{\Pi}_1 \times \check{\Pi}_2, 1-s), $$
where $\epsilon(\Sigma_1 \times \Sigma_2, s)$ is the product of the 
local factors of [JPS].
Moreover, (2.4.1) implies that the archimedean $L$ and $\epsilon$-factors of
$\Sigma_1 \times \Sigma_2$ and $\Pi_1 \times \Pi_2$ are equal.
It then follows from [He1, Theorem 4.1] that, for any non-trivial additive
character $\psi$ of $F$, we have
$$\gamma(\Sigma_{1,v_1} \times \Sigma_{2,v_1},s,\psi) =
\gamma(\Pi_{1,v_1} \times \Pi_{2,v_1},s,\psi). \tag 3.2.3 $$
Here the subscript ${}_{v_1}$ designates the local factor at $v_1$, 
and $\gamma(\bullet,s,\psi)$
is the local ``gamma" factor
$$\gamma(\Pi_{1,v_1} \times \Pi_{2,v_1},s,\psi) =
\frac{\epsilon(\Pi_{1,v_1} \times \Pi_{2,v_1},s,\psi)\cdot 
L(\check{\Pi}_{1,v_1} \times \check{\Pi}_{2,v_1},1-s)}{L(\Pi_{1,v_1} 
\times \Pi_{2,v_1},s)},
\tag 3.2.4 $$
with the analogous formula when the $\Pi$'s are replaced by $\Sigma$'s.

Now it follows from Lemma 3.1 (b) and Theorem 1.7 that 
$\Sigma_{i,v_1} = \sigma(\Pi_{i,v_1})$.
Lemma 3.1(c) implies that $\Sigma_{i,v_1} = \sigma_i \otimes \beta_i$, so that
$\Pi_{i,v_1} = \pi_i \otimes \beta_i$, since the correspondence 
$\sigma(\bullet)$
commutes with twists by characters.  Let $s_0$ be a complex number such that
$\beta_1\cdot \beta_2 = |\bullet|^{s_0}$, where $|\bullet|$ is the 
absolute value character.
Thus we have the identity
$$\gamma(\sigma_1 \times \sigma_2,s+s_0,\psi) = \gamma(\pi_1 \times 
\pi_2,s+s_0,\psi). \tag 3.2.5 $$
Since $\pi_1$ and $\pi_2$ are taken to be supercuspidal, the $L$-factors in
the numerator and denominator of the right-hand side of (3.2.4) have no common
poles (and indeed are both trivial unless $n_1 = n_2$ and $\pi_1 = 
\check{\pi}_2$).
The corresponding fact holds for the left-hand side, because 
$\sigma_1$ and $\sigma_2$
are irreducible.  The equality (3.2.1) then follows as in [He1,\S 4] 
and [LRS, p. 318].
The final assertions are then immediate from the additivity properties of the
local factors and from [He3].
\enddemo

\bigskip
\bigskip

\centerline{\bf 4.  Remarks on the general case}
\bigskip

We have already noted that [BZ] provides an extension of $\pi \mapsto 
\sigma(\pi)$ to
a bijection $\Cal{A}(n,F) \leftrightarrow \Cal{G}(n,F)$, for all $n$.  We again
denote this bijection $\pi \mapsto \sigma(\pi)$.  The inverse bijection
is denoted $\sigma \mapsto \pi(\sigma)$.  If $n = n_1 + \dots + n_r$ is a
partition of $n$, and if $\pi_i \in \CA_0(n_i,F)$, for $1 \leq i \leq r$, then
we write $\pi_1 \boxplus \cdots \boxplus \pi_r \in \CA(n,F)$ 
(Langlands sum) for the
inverse image under $\sigma$ of $\sigma(\pi_1) \oplus \cdots \oplus 
\sigma(\pi_r)$:
$$\pi_1 \boxplus \cdots \boxplus \pi_r = \pi(\sigma(\pi_1)\oplus 
\cdots \oplus \sigma(\pi_r)).$$
Let $\Cal{G}_{ss}(n,F) \subset \Cal{G}(n,F)$ denote the subset of 
representations
of the Weil group (i.e., representations without monodromy operator), and let
$\CA_{ss}(n,F) \subset \CA(n,F)$ denote the corresponding subset.  Then
$\CA_{ss}(n,F)$ can be described as the set of irreducible admissible 
representations
of $GL(n,F)$ obtained as Langlands sums of supercuspidals, and
$\bigcup_n \CA_{ss}(n,F)$ (disjoint union) becomes a monoid under 
Langlands sum.  We let $RG(F)$ denote
the Grothendieck group of representations of $W_F$.  Then $\sigma$ 
places $RG(F)$
in bijection with the group completion $RA(F)$ of the monoid 
$\bigcup_n \CA_{ss}(n,F)$.

The local $\epsilon$-factors $\epsilon(\sigma \times \sigma',s,\psi)$ 
attached to
pairs of Weil group representations are additive in each of the two 
factors $\sigma$,
$\sigma'$, with respect to direct sums [D2, Theorem 4.1].  The same 
is true of the
local $L$-factors, hence of the local gamma-factor of pairs discussed 
in the proof
of Theorem 3.2.  Similarly, and by design, the local $\epsilon$-factors
$\epsilon(\pi \times \pi',s,\psi)$ of [JPS], attached to $\pi \in \Cal{G}(n,F)$
and $\pi' \in \Cal{G}(m,F)$, are additive in $\pi$ and $\pi'$ with respect to
Langlands sums; the same is true of the local gamma-factors.  It thus follows
that $\epsilon(\pi \times \pi',s,\psi)$ (resp. $\gamma(\pi \times 
\pi',s,\psi)$)
extends to a function on $RA(F) \times RA(F)$ with values in the 
multiplicative group of
entire (resp. meromorphic) functions on the complex line; here $\psi$ 
is viewed as fixed.

Let $\sigma \in \Cal{G}(n,F)$, $\sigma' \in \Cal{G}(m,F)$, and define
$$\epsilon_A(\sigma \times \sigma',s,\psi) = \epsilon(\pi(\sigma) 
\times \pi(\sigma'),s,\psi), \tag 4.1$$
with the right-hand side defined as in [JPS].
In this way, the automorphic local $\epsilon$-factor of pairs defines 
a function
on pairs of representations of $W_F$, and indeed on $RG(F) \times 
RG(F)$.  Thus the correspondence
$\pi \mapsto \sigma(\pi)$ qualifies as the conjectured local 
Langlands correspondence
provided
$$\epsilon_A(\sigma \times \sigma',s,\psi) = \epsilon(\sigma \times 
\sigma',s,\psi) \tag 4.2$$
for all pairs $(\sigma,\sigma')$ as above.

With $F$ fixed, we let $G^0(F)$ denote the set of pairs $(F',\chi)$, 
where $F'$ is a finite
extension of $F$ and $\chi$ is a character of $W_{F'}$, or, 
equivalently, of $(F')^{\times}$.
Let $R^0(F)$ denote the free abelian group on the elements of $G^0(F)$.
There is a natural homomorphism
$$\phi: R^0(F) \rightarrow RG(F)$$
defined on generators by
$$\phi((F',\chi)) = Ind_{F'/F}~\chi \in \Cal{G}_{ss}([F':F],F). \tag 4.3$$
By Brauer's theorem, $\phi$ is surjective.  The following lemma is thus
an immediate consequence of our earlier remarks on multiplicativity 
of $\epsilon$-factors.

\proclaim{Lemma 4.4}  Suppose, for every pair 
$((F',\chi_1),(F'',\chi_2)) \in G^0(F) \times G^0(F)$
we have
$$\epsilon_A(Ind_{F'/F}~\chi_1 \times Ind_{F''/F}~\chi_2,s,\psi) = 
\epsilon(Ind_{F'/F}~\chi_1 \times Ind_{F''/F}~\chi_2,s,\psi).$$
Then (4.2) holds for all pairs $(\sigma,\sigma')$.
\endproclaim

Lemma 4.4 admits the following reformulation:

\proclaim{Lemma 4.4.1}  Suppose, for every pair 
$((F',\chi_1),(F'',\chi_2)) \in G^0(F) \times G^0(F)$
we have
$$\gamma(\pi(Ind_{F'/F}~\chi_1) \times 
\pi(Ind_{F''/F}~\chi_2),s,\psi) = \gamma(Ind_{F'/F}~\chi_1 \times 
Ind_{F''/F}~\chi_2,s,\psi), \tag 4.4.2$$
where the left-hand side is defined by (3.2.4) and the right-hand 
side is the Galois-theoretic
analogue.  Then (4.2) holds for all pairs $(\sigma,\sigma')$.
\endproclaim
\demo{Proof}  We have to show that (4.4.2) implies (4.2).
It suffices to show (4.2) when $\sigma$ and $\sigma'$ are irreducible, and thus
$\pi(\sigma)$ and $\pi(\sigma')$ are supercuspidal.  As in the proof 
of Theorem 3.2, (4.2)
then follows from the identity of gamma-factors in the 
irreducible/supercuspidal case.
But (4.4.2) implies the analogue of (4.2) for gamma-factors, by the 
same argument used
to prove Lemma 4.4.
\enddemo

Let $(F',\chi) \in G^0(F)$, with $[F':F] = n$, and let $\tilde{F}$ be the
Galois closure of $F'$ over $F$.  Let $E'/E$ be an extension of
CM fields as in Theorem 1.2 (ii), with $E$ containing the imaginary 
quadratic field $K_0$,
in which $p$ splits, and let $\tilde{E}$ be the
Galois closure of $E'$ over $E$.  We choose complex conjugate places 
$v_1$ and $v_2$ of $E$
dividing $p$, as in \S 2, so that $E_{v_1} \cong F$.  We assume that 
there is exactly
one prime $\tilde{v}_1$ (resp. $\tilde{v}_2$) of $\tilde{E}$ dividing 
$v_1$ (resp. $v_2$),
and that the isomorphism $E_{v_1} \isoarrow  F$ extends
to an isomorphism $\tilde{E}_{\tilde{v}_1} \isoarrow \tilde{F}$.
Let $v_i'$, $i = 1, 2$, denote the corresponding primes of $E'$, so 
that $E'_{v_1'} \cong F'$.

The strategy for proving (4.4.2) is analogous to that used above in 
the tame case.
Our goal is to embed $\chi$ as the local component of an 
$E'/E$-regular Hecke character
$\bold{\chi}$ (note change in notation) for which we can construct a 
weak automorphic
induction $\pi \in \CU(n,E)$.  This is more involved than in the tame case,
since the local Galois groups are not usually so simple, but the 
general principle is the same.
What is missing to complete the argument is the analogue of Theorem 1.7.  More
precisely, it will generally not be the case that the local component 
$\pi_{v_1}$
is supercuspidal.  Thus we have no information about the relation between
$\sigma(\pi_{v_1})$, defined by extension of the supercuspidal/irreducible
correspondence, and $Ind_{F'/F}~\chi$.

Since $Ind_{F'/F}~\chi$ is no longer assumed to be irreducible, our techniques
have to be modified.  We assume there is a second pair $v(*)_1$ and $v(*)_2$ of
complex conjugate places of $E$ dividing $p$, and that each $v(*)_i$, 
$i = 1, 2$, is divisible by
exactly one prime $\tilde{v}(*)_i$ of $\tilde{E}$, with compatible
isomorphisms
$$E_{v(*)_1} \isoarrow  F, \qquad \tilde{E}_{\tilde{v}(*)_1} 
\isoarrow \tilde{F}.$$
This can be arranged by replacing $E$ by its compositum
with a real quadratic field in which $p$ splits.
Define $v(*)_i'$, $i = 1, 2$, as before, so that $E'_{v(*)_1'} \cong F'$.

\proclaim{Definition 4.5}  Let $\chi(*)$ be a complex-valued 
character of $W_{F'}$, or equivalently a
character of $F^{\prime,\times}$.  We say that $\chi(*)$ is {\bf in 
general position}
(relative to $F$) if there
is a sequence $\tilde{F} = F_0 \supset F_1 \supset \dots \supset 
F_{t-1} \supset F_t = F$
such that

(i) For all $i$, $F_i/F_{i+1}$ is a cyclic extension of prime degree;

(ii) The Mackey constituents of
$Ind_{F'/F} ~\chi(*)$ restricted to $W_{F_i}$ are all irreducible, 
for $i = 0, \dots, t$.
\endproclaim

The Mackey constituents have been defined in the Notation section.  For
$a \in A_i = W_{F'} \backslash W_F/W_{F_i}$ we let
$I(a,i,\chi(*)) = Ind_{a(F')\cdot F_i/F_i} ~a(\chi(*))$ denote the
corresponding Mackey constituent, where $a(F')$ is viewed as a subfield
of $\tilde{F}$.  Let $n(a,i) = [F_i:F_i\cap a(F)]$ and let
$\pi(a,i,\chi(*)) \in \CA_0(n(a,i),F_i)$ denote the supercuspidal 
representation
that corresponds to $I(a,i,\chi(*))$.

Let $\sigma(*) = Ind_{F'/F} ~\chi(*) = I(1,t,\chi(*))$.  If $\chi(*)$ 
is in general
position then $\sigma(*)$ is irreducible and corresponds to a 
supercuspidal representation
$\pi(*) = \pi(\sigma(*))$ of $GL(n,F)$.  Moreover, for each $i$ we have
$$BC_{F_{i}/F_{i+1}} \circ BC_{F_{i+1}/F_{i+2}} \circ \cdots \circ 
BC_{F_{t-1}/F}~\pi(*) \cong
\boxplus_{a \in A_i} \pi(a,i,\chi(*)). \tag 4.6$$

The proof of the following lemma, which is probably well known
to experts, is postponed to the end of this section.

\proclaim{Lemma 4.7}  Let $F'/F$ be an extension of local fields.  There exist
characters of $F^{\prime,\times}$ of arbitrarily large conductor that are in
general position relative to $F$.
\endproclaim

Now we can prove the following extension of Proposition 2.4:

\proclaim{Proposition 4.8}  Let $\chi$ and $\chi(*)$ be characters of $W_{F'}$,
with $\chi(*)$ in general position relative to $F$.  Let $\bc$ be an 
$E'/E$-regular
algebraic Hecke character of $E'$ with components $\chi$ and $\chi(*)$ at
$v'_1$ and $v(*)_1'$, respectively, and such that $\bc = \bc^{-1} 
\circ c$.  Then
there is a cuspidal
automorphic representation $\Pi \in \CU(n,E)$ which is a weak 
automorphic induction of $\bc$.
Moreover, the local component $\Pi_{v(*)_1}$ at $v(*)_1$ is supercuspidal.
\endproclaim

\demo{Proof} As before, when referring to representations of global 
Weil groups we
use the language of complex representations and their associated $\lambda$-adic
families interchangeably.

The isomorphism $Gal(\tilde{E}/E) \isoarrow Gal(\tilde{F}/F)$
provides a collection of subfields $E_i \subset \tilde{E}$ with
$Gal(\tilde{E}/E_i) = Gal(\tilde{F}/F_i)$, $i = 0, \dots, t$.
We prove for each $i$ the existence of an isobaric representation
$\Pi_i \in \CU(n,E_i)$ whose weakly associated $\lambda$-adic family is given
by the restriction to $Gal(\overline{E}_i/E_i)$ of $Ind_{E'/E} ~\bc$.
Note that $A_i$ gives a set of representatives of
$Gal(\overline{E}/E')\backslash Gal(\overline{E}/E) /Gal(\overline{E}/E_i)$.
Thus the Mackey constituents of  $Ind_{E'/E} ~\bc$, restricted to
$Gal(\overline{E}_i/E_i)$, are parametrized by $A_i$.  For $a \in A_i$ we
let  $I(a,i,\bc)$ denote the corresponding Mackey constituent.

I claim that, by induction on
$t$, I may assume that $I(a,i,\bc)$ is weakly associated to a cuspidal
automorphic representation $\Pi(a,i) \in \CU(n(a,i),E_i)$, for $i = 
0, \dots, t-1$
and all $a \in A_i$, with local component supercuspidal
at the prime dividing $v(*)_1$.  Indeed, this is clear for $i = 0$, 
since in this case
$n(a,i) = 1$ for all $a \in A_i$.  Moreover, the degree
$[E_0:E_i]$ is a proper divisor of $[E_0:E]$ for $i < t$. To prove 
the claim, it thus
suffices to show that, for all $a \in A_i$, the Hecke character
$a(\bc) \circ N_{a(E')\cdot E_i/a(E')}$ of $a(E')\cdot E_i$
satisfies the same hypotheses relative to $E_i$ as did $\bc$ relative 
to $E$, namely:

(i)  $a(\bc)\circ N_{a(E')\cdot E_i/a(E')}$ is $a(E')\cdot E_i /E_i$-regular;

(ii) $a(\bc)\circ N_{a(E')\cdot E_i/a(E')} = [a(\bc)\circ 
N_{a(E')\cdot E_i/a(E')}]^{-1} \circ c$;

(iii) The local constituent of $a(\bc)\circ N_{a(E')\cdot E_i/a(E')}$ 
at the prime
dividing $v(*)_1$ is in general position relative to $F_i$.

Of these hypotheses, (i) and (ii) are obvious, and (iii) follows from 
the hypothesis
of general position for $\chi(*)$, since the Mackey constituents of 
the restriction of
$I(a,i,\bc)$ to $Gal(\overline{E}_j/E_j)$, for $j < i$, are among the Mackey
constituents of $Ind_{E'/E} ~\bc$, restricted to $Gal(\overline{E}_j/E_j)$
(``transitivity of restriction").

\noindent {\bf Remark.}  The preceding argument makes use of the full 
strength of
Theorem 1.7, in that the infinite component of $\Pi(a,i)$ is a priori 
an arbitrary
representation of cohomological type.  For any archimedean place $\tau$ of
$E$ it is possible to choose the infinity type of $\bc$ so that, for 
a specific choice
of place $\tau_i$ dividing $\tau$, each $\Pi(a,i)_{\tau_i}$ is a cohomological
representation with the infinitesimal character of a one-dimensional 
representation; i.e.,
an abelian twist of the Langlands sum considered in [H1,(22)]. 
However, I see no way
to control the local components of $\Pi(a,i)$ at the other primes of $E_i$
dividing $\tau$.

Thus by induction we may define
$$\Pi_{t-1} = \boxplus_{a \in A_{t-1}} \Pi(a,t-1) \in \CU(n,E_{t-1}),$$
whose weakly associated $\lambda$-adic family $\Sigma_{t-1}$ is given
by the restriction to \linebreak $Gal(\overline{E}_{t-1}/E_{t-1})$ of 
$Ind_{E'/E} ~\bc$.
Now $E_{t-1}/E$ is cyclic of prime degree $q$, say, and $Gal(\overline{E}/E)$
acts transitively on the right on $A_{t-1}$ 
($Gal(\overline{E}/E_{t-1})$ is a normal subgroup
of $Gal(\overline{E}/E)$).  Moreover, $\Sigma_{t-1}$ is the 
restriction of a $\lambda$-adic family
of representations of $Gal(\overline{E}/E)$, hence is 
$Gal(E_{t-1}/E)$-invariant.  The
weakly associated isobaric representation $\Pi_{t-1}$ is thus also 
$Gal(E_{t-1}/E)$-invariant,
by strong multiplicity one.

Thus there are two
possibilities.  If $A_{t-1}$ has $q$ elements then it is a principal 
homogeneous
space under $Gal(E_{t-1}/E)$.  Choosing a basepoint $e \in A_{t-1}$ 
and denoting
by $\alpha$ a generator of $Gal(E_{t-1}/E)$, we thus have
$$\Pi_{t-1} = \boxplus_{a \in A_{t-1}} \Pi(a,t-1) =
\Pi(e,t-1)\boxplus \alpha(\Pi(e,t-1)) \boxplus \cdots \boxplus 
\alpha^{q-1}(\Pi(e,t-1)).$$
Then $\Pi = AI_{E_{t-1}/E} ~\Pi(e,t-1)$ is the unique automorphic 
representation
of $GL(n,E)$ that base changes to $\Pi_{t-1}$ [AC,Lemma III.6.4]. 
Moreover, it follows
from the irreducibility of $Ind_{F'/F}~\chi(*)$ that the 
supercuspidal representations
$\alpha^j(\Pi(e,t-1)_{v(*)_1})$ are all distinct, hence 
[HH,Proposition 5.5] that the local
component $\Pi_{v(*)_1}$ is supercuspidal.  It is thus clear that $\Pi$ is a
weak automorphic induction of $\bc$.

It remains to consider the case in which $A_{t-1}$ is a singleton 
$e$.  Then $\Pi_{t-1} = \Pi(e,t-1)$
is already cuspidal and has $q$ distinct descents to cuspidal automorphic
representations of $GL(n,E)$.  On the other hand, $I(e,t-1,\bc)$ has 
$q$ distinct extensions
to irreducible $\lambda$-adic families of representations of 
$Gal(\overline{E}/E)$.  Thus
we conclude by applying the argument used to prove Lemma 1.6.
\enddemo

To continue, we need to state a conjecture.

\proclaim{Conjecture 4.9}  Let $E$ be a CM field as in Theorem 1.2 (ii) and let
$\Pi \in \CU(n,E)$.  Let $\Sigma$ be the $\lambda$-adic family weakly 
associated
to $\Pi$, and let $w$ be a finite place of $E$.  Let $\Sigma_w$ denote the
restriction of $\Sigma$ to a decomposition group at $w$.  Then
$$\Sigma_w \cong \sigma(\Pi_w).$$
\endproclaim

This conjecture was proved by Carayol when $n = 2$ [Ca], with
$\sigma(\Pi_w)$ given by the local Langlands correspondence.  In the general
case, it can be translated
into a problem about bad reduction of the Shimura varieties 
considered in [C2] and
[H1].  Translation into a precise problem is probably the main step in proving
the conjecture.  Theorem 1.2 states that the conjecture is true for almost all
unramified places $w$, and Theorem 1.7 states that it is true at 
supercuspidal places.
A weaker version, probably much easier to prove, is sufficient for 
our purposes:

\proclaim{Conjecture 4.9 bis}  Under the hypotheses of Conjecture 4.9,
let $\Sigma_{w,ss} \in \Cal{G}_{ss}(n,F)$ denote the restriction of $\Sigma_w$
to $W_F$ (i.e., forget about monodromy), and define 
$\sigma(\Pi_w)_{ss}$ analogously.  Then
$$\Sigma_{w,ss} \cong \sigma(\Pi_w)_{ss}.$$
\endproclaim

\proclaim{Proposition 4.10}  Suppose Conjecture 4.9 bis.  Then the 
identity of $\epsilon$-factors
(4.2) holds for all pairs $(\sigma,\sigma')$.
\endproclaim

In other words, the local Langlands conjecture would follow from 
Conjecture 4.9 bis.  We will
see from the proof that it suffices to know Conjecture 4.9 bis when 
$\Sigma_{w}$ is a
{\it monomial} representation.

\demo{Proof}  It follows from Lemma 4.4.1 that it suffices to prove 
the identity of local
gamma-factors (4.4.2) for pairs of monomial representations 
$(\phi(F_1,\chi_1),\phi(F_2,\chi_2))$.
Now consider the analogue of Lemma 3.1 in our situation:  there 
exists a global CM field
$E = E_0\cdot K_0$ as in Lemma 3.1, CM extensions $E_1$ and $E_2$ of 
$E$, $E_i/E$-regular
Hecke characters $\bc_i$, and cuspidal automorphic representations
$\Pi_i \in \CU(n_i,E)$ such that $\Pi_i$ is a weak automorphic 
induction of $\chi_i$,
$i = 1, 2$, satisfying the hypothesis of Lemma 3.1 (c).  This 
generalization of Lemma 3.1
follows from Proposition 4.8, by the argument preceding Lemma 3.1. 
Then the argument
used to prove Theorem 3.2 yields the identity
$$\gamma(\Pi_{1,v_1} \times \Pi_{2,v_1},s,\psi) = 
\gamma(\phi(F_1,\chi_1) \times \phi(F_2,\chi_2),s,\psi), \tag 4.11$$
By Conjecture 4.9 bis we can replace $\Pi_{1,v_1} \times \Pi_{2,v_1}$ by
$\pi(\phi(F_1,\chi_1)) \times \pi(\phi(F_2,\chi_2))$ on the left-hand 
side.  This has
the effect of transforming (4.11) into (4.4.2) and thus completes the proof.
\enddemo

It remains to prove Lemma 4.7.  We will construct characters 
$\chi(*)$ of $F^{\prime,\times}$
of finite order that are in general position relative to $F$, with 
arbitrarily large
conductor.  In particular, in the Galois-theoretic notation of \S 2, we
write $G = \Gamma_F$, $H = \Gamma_{F'}$, $\tilde{G} = 
\Gamma_{\tilde{F}}$; $\chi(*)$ will be viewed equivalently
as a character of $F^{\prime,\times}$ or of $H$.  We first describe a 
sufficient
condition for $Ind_H^G ~\chi(*)$ to be irreducible.  Let $A \subset 
G$ be a set of representatives
for $H\backslash G/H$, with $e$ the representative for the identity double
coset.  By Clifford-Mackey theory, for $Ind_H^G ~\chi(*)$ to be irreducible
it is necessary and sufficient that, for all $a \in A$, $a \neq e$, 
$\chi(*)$ and
$a(\chi(*))$ have distinct restrictions to $H\cap aHa^{-1}$.  In particular,
letting $\tc = \chi(*) \circ N_{\tilde{F}/F'}$, a sufficient condition for
irreducibility of $Ind_H^G ~\chi(*)$ is that $H$ is the stabilizer in 
$G$ of $\tc$:
$$H = \{g \in G| g(\tc) = \tc\}. \tag 4.12$$

On the other hand, let $K$ be any subgroup of $G$ containing $\tilde{G}$.
Recall the explicit description of the Mackey constituents of 
$Ind_H^G ~\chi(*)$,
restricted to $K$:
$$Ind_{aHa^{-1}\cap K}^K ~a(\chi(*)),~~ a \in H\backslash G/K.$$
It follows easily from Clifford-Mackey theory that (4.12) is a 
sufficient condition
for $\chi(*)$ to be in general position relative to $F$.

So it suffices to construct characters $\chi(*)$ with arbitrarily 
large conductor
satisfying (4.12).  Let $X(F')$, resp. $X(\tilde{F})$ denote the 
groups of characters
of $F^{\prime,\times}$, resp. of $\tilde{F}^{\times}$, with values in 
$\CC^{\times}$, and let
$\nu:  X(F') \rightarrow X(\tilde{F})$ denote pullback via $N_{\tilde{F}/F'}$.
Similarly, let $X_p(F')$ and $X_p(\tilde{F})$ denote the groups of 
$\ZZ_p^{\times}$-valued
characters, and $\nu_p: X_p(F') \rightarrow X_p(\tilde{F})$ the pullback.
It follows from class field theory that the image of $\nu$ is of finite
index in $X(\tilde{F})^H$.  Thus it suffices to show that $X(\tilde{F})^H$
contains characters $\xi$ of arbitrarily large conductor such that
$$H = \{g \in G| g(\xi) = \xi\}. \tag 4.13$$

In fact, it is enough to find characters $\xi_p$ in 
$X_p(\tilde{F})^H$ of infinite
conductor satisfying (4.13).  Indeed, if $\xi_p$ satisfies (4.13), then so does
its reduction modulo $p^N$:
$$\xi_p \pmod{p^N}: \tilde{F}^{\times} \rightarrow (\ZZ/p^N \ZZ)^\times$$
for sufficiently large $N$. By further increasing $N$, we can guarantee that
the conductor of $\xi_p \pmod{p^N}$ is arbitrarily large;
then composing with an embedding $(\ZZ/p^N \ZZ)^\times 
\hookrightarrow \CC^{\times}$,
we obtain a complex character with the desired property.

Let $M \subset \tilde{F}$ denote the image of $\tilde{F}^{\times}$ under
the $p$-adic logarithm map, and let $\Cal{O}_F$ denote the ring of 
integers in $F$.
For any subgroup $M' \subset M$ we let $X_p(M')$ denote the
group of $\ZZ_p^{\times}$-valued characters of $M'$, and
$X_p^1(M') \subset X_p(M')$ the subgroup of characters with values in
$1+p\ZZ_p$.  It suffices to find a
$G$-invariant subgroup $M'$ of $M$ of finite index and a character 
$\xi \in X_p^1(M')$
satisfying (4.13).  Indeed, such a character necessarily has infinite 
conductor.
Multiplication by $p^h$ for sufficiently large $h$ defines a $G$-equivariant
embedding of $M$ in $M'$, hence a restriction $r_h: X_p^1(M')^H 
\rightarrow X_p^1(M)^H$.
The $X_p^1$ groups are free $\ZZ_p$-modules.  Hence $r_h$ is injective and
$r_h(\xi)$ again satisfies (4.13).  The same is therefore true
of its inflation to $X_p(\tilde{F})^H$.

But now, by using the normal basis theorem for
the Galois extension $\tilde{F}/F$, we see easily that $M$ contains
a subgroup $M'$ of finite index, isomorphic as $\ZZ_p[G/\tilde{G}]$-module to
$\Cal{O}_F[G/\tilde{G}]$. By duality,
$$X_p^1(M') \simeq Hom(\Cal{O}_F[G/\tilde{G}],\ZZ_p)$$
also contains a subgroup isomorphic as $\ZZ_p[G/\tilde{G}]$-module to 
$\Cal{O}_F[G/\tilde{G}]$.
The condition (4.13) is dense in $\Cal{O}_F[G/\tilde{G}]^H$, being the
complement of a finite union of proper $\Cal{O}_F$-submodules.  Thus
there are $p$-adic characters of $M'$ of infinite conductor satisfying (4.13).
This completes the proof of Lemma 4.7.

\bigskip

\centerline{\bf  Appendix. Remarks on Theorem 1.7.}

As promised at the end of \S 1, we sketch how to extend the results
of [H1] on the compatibility of the local Galois correspondence
with the global correspondence realized on the cohomology of certain
Shimura varieties, as needed for the applications of the present 
paper.  We freely
make use of the notation and techniques of [H1] and [RZ].

In [H1] we work with a Shimura variety denoted $S(G\Cal{G},X_{n-1})_{C(N)}$
attached to the $\QQ$-group $G\Cal{G} = GU(\Cal{B},\alpha)$, the 
unitary similitude group of a division
algebra $\Cal{B}$ over the CM field $E$ with involution $\alpha$ of 
the second kind.
Here $C(N)$ is a level subgroup, depending on a positive integer $N$, and
assumed to factor as the product
$\prod_q C(N)_q$ over the rational primes $q$, and the similitude factor is
assumed to be rational.
The field $E$, denoted $\Cal{K}$ in [H1], is assumed to be of the 
form  $E_0 \cdot K_0$
as in Theorem 1.2 (ii), and $p$ is assumed to split in $K_0$ as the product
$\frak{p}_1 \cdot \frak{p}_2$.  The primes of $E$ dividing 
$\frak{p}_i$ are denoted
$v_i^{(j)}$, $i = 1, 2$, $j = 1, \dots, s$, with $v = v_1^{(1)}$ the
distinguished place for which $E_v = F$.  In general, we let
$F_j$ denote the completion of $E$ at $v_i^{(j)}$.  Then the $p$-adic 
points of $G\Cal{G}$
are given by
$$G\Cal{G}(\QQ_p) \cong \prod_j \Cal{B}(F_j)^{\times} \times 
\QQ_p^{\times} \tag A.1 $$
(cf. [H1,(2)]).

Starting in \S 2 of [H1]  $\Cal{B}$ was
assumed to be a division algebra at every $v_i^{(j)}$.  Moreover,
the $p$-level subgroup $C(N)_p$ was assumed to have a factorization
$$C(N)_p \cong C_v^{N,0} \times \prod_{j > 1} C_p^{(j)} \times 
\ZZ_p^{\times}, \tag A.2 $$
where $C_v^{N,0}$ is a principal congruence subgroup defined in [loc. cit.] but
$C_p^{(j)}$ is assumed {\it maximal} for $j > 1$.
These assumptions were
made exclusively to simplify notation, and to be able to refer freely 
to the discussion
in Chapter 6 of [RZ].  Upon closer inspection, [RZ] turns out to assume that
$\Cal{B}$ is {\it split} at every $v_i^{(j)}$ for $j > 1$, but the 
corresponding
$C_p^{(j)}$ are still assumed maximal.  In any case, this restriction 
is irrelevant,
and the general case can be found in the literature:  implicitly in the earlier
chapters of [RZ] and explicitly in [BoZ] and [V].  The non-expert 
will be bewildered to
discover that no two of these three references work with quite the 
same Shimura datum
$(G\Cal{G},X_{n-1})$.  In [H1] and [RZ] the Shimura variety parametrizes
weight $-1$ Hodge structures; in [BoZ] and [V] the weights are $1$ and $0$,
respectively.  Moreover, [BoZ] and [V] both use the full similitude 
group rather than
the group with rational similitude factor.  Passage between these 
points of view
is standard for specialists in Shimura varieties and we will say no more on this 
point.

Now suppose we are given a supercuspidal representation $\pi_v$ of $GL(n,F)$
and an automorphic representation $\pi$ of $G\Cal{G}$ of 
cohomological type whose $p$-adic
constituent $\pi_p$ factors with respect to (A.1) as
$$\pi_v \otimes (\bigotimes_{j > 1} \pi_j) \otimes \eta_p.  \tag A.3 $$
The reader will verify that, as in [H1, p. 95 ff.], we may assume the 
character $\eta_p$
of $\QQ_p^{\times}$ to be trivial.   However, in \S 4 above we allowed at least
one $\pi_j$ to be non-trivial, and indeed quite general.   It needs 
to be established that
the local Galois representation $\sigma(\pi_v)$ attached to $\pi$ by the recipe
of [H1,\S3] depends only on $\pi_v$ and not on the other local 
factors of $\pi$.
The argument from $p$-adic uniformization (as in [H1, Proposition 2] 
and recalled below)
shows easily that $\sigma(\pi_v)$ is independent of $\pi_q$ for $q \neq p$.
Thus it is a matter of showing that $\sigma(\pi_v)$ depends neither 
on the invariants
of $\Cal{B}(F_j)$ nor on $\pi_j$ for $j > 1$.

The definition of $\sigma(\pi_v)$ has been recalled in the Erratum to 
[H1].  Replacing
the notation $G\pi$ used there by the current notation $\pi$, we have
$$ \sigma(\pi_v) = [\tilde{\sigma}(\pi_v) \otimes 
\bar{\nu}(\pi)_v^{-1}]^*, \tag A.4$$
where $*$ denotes contragredient and where
$$
\tilde{\sigma}(\pi_v) =
[Hom_{GJ}(H^{n-1}_c(\cbM_N,\Qbar_{\ell})_{SS(F)}, \pi_p) \otimes 
GJL(\pi_v)^{*}]^{GG}.  \tag A.5$$
An elementary calculation
shows that $\bar{\nu}(\pi)_v$ may well depend on the character 
$\eta_p$ (which we have
assumed trivial) but is independent of the $\pi_j$ for $j > 1$.  Moreover,
the factor $GJL(\pi_v)^{*}$ in (A.5) does not affect the Galois 
action.  It thus remains
to show that the dependence of the $\GalEv$ module
$$Hom_{GJ}(H^{n-1}_c(\cbM_N,\Qbar_{\ell}), \pi_p)\otimes \nu(\pi)_v \tag A.6$$
depends only on $\pi_v$.  Here $GJ$ is the group of $p$-adic points 
of the inner
twist $G\Cal{J}$ of $G\Cal{G}$ used in [H1]:
$$GJ = GL(n,F) \times \prod_{j >1} \Cal{B}(F_j)^{\times} \times 
\QQ_p^{\times}. \tag A.7$$

The rigid parameter space $\cbM_N$ actually depends on the level 
subgroup $C_p$, which
in the present paper we are allowing to vary.  More precisely, we take $\cbM_N$
to be the inverse limit as $\prod_{j > 1} C_p^{(j)}$ shrinks to the 
trivial group, while $C_v^{N,0}$
remains fixed.  Just as in the appendix to [H1], $\cbM_N$ is 
contained in a product
$\prod_j \cbM_j$.  Here $\cbM_1 = \cbM_{v,N}$, is Drinfeld's rigid 
space of level $N$
attached to $F$, normalized as in [H1] (following [RZ])
to include a morphism to $\QQ_p^{\times}/\ZZ_p^{\times}$ for the 
similitude factor (polarization).
This morphism splits (non-canonically) to yield an isomorphism
$$\cbM_1 \isoarrow \cbM_{v,N}^+ \times \QQ_p^{\times}/\ZZ_p^{\times}.$$
For $j > 1$ there is a non-canonical identification
$\cbM_j \isoarrow \Cal{B}(F_j)^{\times} \times 
\QQ_p^{\times}/\ZZ_p^{\times}$, and
$\cbM_N \subset \prod_j \cbM_j$ can be defined as the subset of the product on
which the natural maps to $\QQ_p^{\times}/\ZZ_p^{\times}$ agree. 
Thus we can identify
$$\aligned
  \cbM_N &\isoarrow \prod_j \cbM_{v,N}^+ \times \prod_{j > 1} 
\Cal{B}(F_j)^{\times}
\times \QQ_p^{\times}/\ZZ_p^{\times}.\\
&\isoarrow ~~~\cbM_{v,N} \times  \prod_{j > 1} \Cal{B}(F_j)^{\times}
\endaligned \tag A.8$$
The first factorization is compatible with the action of $GJ$ via the 
factorization (A.7);
the second factorization groups together the first and last factors 
of the first.
All this is proved just as in the last chapter of [RZ], or as in [BoZ, \S 1].

Now we recall that $\cbM_N$ is naturally defined over the field $F^{nr}$, the
maximal unramified extension of $F$, but that it is given with a Weil descent
datum [RZ, p. 99]
$$\zeta:  \cbM_N \rightarrow \phi^* \cbM_N$$
([RZ, Proposition 6.49]; [BoZ, pp. 33-34]).
Here $\phi$ is induced from Frobenius acting on $F^{nr}$.  Now 
$Gal(\overline{F}/F^{nr})$
acts trivially on the last two factors of (A.8).  Moreover, the 
explicit calculation
of the Weil descent data in [RZ] and [BoZ] shows that, under the second factorization
of (A.8), the last term splits off:
$$\zeta = \zeta_v \times 1: \cbM_{v,N} \times \prod_{j > 1} 
\Cal{B}(F_j)^{\times}
\rightarrow \phi^*(\cbM_{v,N}) \times  \prod_{j > 1} 
\Cal{B}(F_j)^{\times}. \tag A.9$$
In other words, the factor $\prod_{j > 1} \Cal{B}(F_j)^{\times}$ is irrelevant
to the Galois representation on $H^{n-1}_c(\cbM_N,\Qbar_{\ell})$.  Thus
$\sigma(\pi_v)$ really does depend only on $\pi_v$.

The calculations that lead to (A.9) are based on the following
considerations.  Suppose for definiteness that
$\Cal{B}(F_j) \cong GL(n,F_j)$; the more general case is analogous.
The moduli space $\cbM_j$, for $j > 1$, parametrizes pairs consisting 
of $\QQ_p$-homogeneous
polarizations and compatible $p$-adic level structures on the $p$-divisible
$O_{F_j}$-module $\bold{X} \times \hat{\bold{X}}$, where $\bold{X}$ 
is the trivial
\'etale $p$-divisible $O_{F_j}$-module of rank $n$ and $\hat{}$ denotes Cartier
dual.  Since the polarization is assumed to be compatible with the 
level structure,
this boils down to a pair consisting of an $O_{F_j}$-linear level 
structure on $\bold{X}$ and a
quasi-isogeny from $\Gm$ to itself.  The Galois group 
$Gal(\overline{F}_j/F_j)$ obviously acts
trivially on the level structure on $\bold{X}$ and by a character on 
the quasi-isogeny.
(The descent datum only becomes effective
when the factor $\QQ_p^{\times}/\ZZ_p^{\times}$ in (A.8) is replaced 
by a finite quotient,
cf. [H1], pp. 114-115.)  This translates directly into (A.9).

\medskip

\centerline{\bf  REFERENCES}
\bigskip
\smallskip

\noindent [AC] Arthur, J. and L. Clozel:  Simple Algebras, Base Change, and the
Advanced Theory of the Trace Formula, {\it Annals of Math Studies} {\bf 120}
(1989).

\smallskip

\noindent [BZ]  Bernstein, J. and A.V. Zelevinski:  Representations 
of the group
$GL(n,F)$, where $F$ is a non-archimedean local field, {\it Russian 
Math. Surveys}
{\bf 31} 1-68 (1976).
\smallskip
\noindent [B]  Blasius, D: Automorphic forms and Galois representations:  some
examples, in L. Clozel and J. S. Milne (eds.), {\it Automorphic 
Forms, Shimura Varieties,
and L-functions},
{\it Perspectives in Mathematics}, {\bf 10}, Vol. II, 1-13 (1990).

\smallskip
\noindent [BR]  Blasius, D. and J. Rogawski:  Motives for Hilbert 
modular forms,
{\it Invent. Math.}, {\bf 114}, 55-87 (1993).

\smallskip

\noindent [BoZ] Boutot, J.-F. and T. Zink:  The p-adic unformization of
Shimura curves, Preprint 95-107, Universit\"at Bielefeld (1995).
\smallskip

\noindent [BF]  Bushnell, C. and A. Fr\"ohlich:  Gauss sums and 
$p$-adic division
algebras, {\it Lecture Notes in Math.} {\bf 987} (1983).
\smallskip

\noindent [BHK]  Bushnell, C., G. Henniart, and P. Kutzko: 
Correspondance de Langlands
locale pour $GL_n$ et conducteurs de paires, manuscript (1997).

\smallskip

\noindent [Ca]  Carayol, H.: Sur les repr\'esentations $\ell$-adiques 
associ\'ees
aux formes modulaires de Hilbert, {\it Ann. scient. Ec. Norm. Sup}, 
{\bf 19}, 409-468 (1986).
\smallskip

\noindent [C1]  Clozel, L.: Motifs et formes automorphes: 
applications du principe
de fonctorialit\'e, in L. Clozel and J. S. Milne (eds.)
{\it Automorphic Forms, Shimura Varieties, and L-functions},
{\it Perspectives in Mathematics}, {\bf 10}, Vol. I, 77-160 (1990).
\smallskip

\noindent [C2]  Clozel, L.: Repr\'esentations Galoisiennes 
associ\'ees aux repr\'esentations
automorphes autoduales de GL(n), {\it Publ. Math. I.H.E.S.}, {\bf 
73}, 97-145 (1991).

\smallskip
\noindent [C3]  Clozel, L.: On the cohomology of Kottwitz's arithmetic
varieties, {\it Duke Math. J.}, {\bf 72}, 757-795 (1993).

\smallskip
\noindent [CH]  Corwin, L. and R. Howe:  Computing characters of 
tamely ramified
$p$-adic division algebras, {\it Pacific J. Math.}, {\bf 73}, 461-477 (1977).

\smallskip
\noindent [D1]  Deligne, P.:  Formes modulaires et repr\'esentations de
$GL(2)$, in P. Deligne and W. Kuyk (eds.), Modular functions of one 
variable, II,
{\it Lecture Notes in Math.} {\bf 349}, 55-105 (1973).

\smallskip
\noindent [D2]  Deligne, P.:  Les constantes des \'equations 
fonctionnelles des fonctions
$L$, in P. Deligne and W. Kuyk (eds.), Modular functions of one variable, II,
{\it Lecture Notes in Math.} {\bf 349}, 501-597 (1973).

\smallskip
\noindent [H1]  Harris, M.: Supercuspidal representations in the cohomology of
Drinfel'd upper half spaces; elaboration of Carayol's program, {\it 
Invent. Math.},
{\bf 129}, 75-119 (1997).

\noindent [H2]  Harris, M.: Galois properties of automorphic representations
of $GL(n)$, in preparation.
\smallskip

\noindent [H3]  Harris, M.: L-functions and periods of polarized 
regular motives,
{\it J. Reine Angew. Math.}, {\bf 483}, 75-161 (1997).
\smallskip

\noindent [He1]  Henniart, G.:  On the local Langlands conjecture for 
$GL(n)$: the
cyclic case, {\it Ann. of Math.}, {\bf 123}, 145-203 (1986).

\smallskip

\noindent [He2]  Henniart, G.:  La conjecture de Langlands locale num\'erique
pour $GL(n)$, {\it Ann. scient. Ec. Norm. Sup}, {\bf 21}, 497-544 (1988).
\smallskip
\noindent [He3]  Henniart, G.: Caract\'erisation de la correspondence 
de Langlands locale
par les facteurs $\epsilon$ de paires, {\it Invent. Math.}, {\bf 113}, 339-350
(1993).
\smallskip

\noindent [He4]  Henniart, G.: Letter, January 1994.

\smallskip
\noindent [HH]  Henniart, G. and R. Herb:  Automorphic induction for 
$GL(n)$ (over
local non-archimedean fields), to appear.

\smallskip
\noindent [JPS]  Jacquet, H., I. I. Piatetski-Shapiro, and J. 
Shalika:  Rankin-Selberg
convolutions, {\it Am. J. Math.}, {\bf 105}, 367-483 (1983).

\smallskip
\noindent [K]  Kazhdan, D.:  On lifting, in {\it Lie Group Representations},
{\it Lecture Notes in Math.}, {\bf 1041}, 209-249 (1984).

\smallskip
\noindent [KZ]  Koch, H. and E.-W. Zink:  Zur Korrespondenz von Darstellungen
der Galoisgruppen und der zentralen Divisionsalgebren \"uber lokalen K\"orpern
(Der zahme Fall), {\it Math. Nachr.}, {\bf 98}, 83-119 (1980).

\smallskip

\noindent [L1]  Langlands, R.P.:  On the functional equation of 
Artin's L-functions,
unpublished manuscript.

\smallskip
\noindent [L2]  Langlands, R.P.:  Automorphic representations, 
Shimura varieties, and
motives.  Ein M\"archen, in {\it Automorphic forms, representations, 
and $L$-functions},
{\it Proc. Symp. Pure Math. AMS}, {\bf 33}, part 2, 205-246 (1979).

\smallskip
\noindent [LRS]  Laumon, G., M. Rapoport, and U. Stuhler, $\Cal{D}$-elliptic
sheaves and the Langlands correspondence, {\it Invent. Math.}, {\bf 
113}, 217-338
(1993).
\smallskip

\noindent [M] Moy, A.:  Local constants and the tame Langlands correspondence,
{\it Am. J. Math.}, {\bf 108}, 863-930 (1986).
\smallskip

\noindent [RZ] Rapoport, M. and T. Zink: Period spaces for 
$p$-divisible groups,
{\it Annals of Math. Studies}, {\bf 141}  (1996).
\smallskip

\noindent [T] Taylor, R.:  Letter to L. Clozel, December 1991.

\noindent [V]  Varshavsky, Y.  $P$-adic uniformization of unitary 
Shimura varieties,
Manuscript (1995).

\smallskip

\baselineskip=9pt
\eightpoint

\hskip 3.6 in U.F.R. de Math\'ematiques

\hskip 3.6 in U.R.A. 748 du C.N.R.S.

\hskip 3.6 in Universit\'e Paris 7

\hskip 3.6 in 2, Pl. Jussieu

\hskip 3.6 in 75251 Paris cedex 05

\hskip 3.6 in FRANCE

\end

\end